\documentclass[a4paper,12pt]{article}
\usepackage{latexsym,amssymb,amsmath,amsthm,amsxtra,xypic}
\usepackage{amsfonts}
\usepackage{amscd}
\usepackage{mathrsfs}
\usepackage[dvips]{graphicx}
\usepackage[all]{xy}

\newcommand{\allowpagebreak}
\allowdisplaybreaks

\newtheorem{thm}{Theorem}[section]
\newtheorem{lem}[thm]{Lemma}

\newtheorem{pro}[thm]{Proposition}
\newtheorem{ex}[thm]{Example}

\theoremstyle{definition}
\newtheorem{definition}[thm]{Definition}

\newcommand {\emptycomment}[1]{}

\newcommand{\pf}{\noindent{\bf Proof.}\ }

\newcommand{\lam}{\lambda}
\newcommand{\si}{\sigma}

\newcommand{\trl}{\triangleleft}
\newcommand{\trr}{\triangleright}
\newcommand{\ppl}{\leftharpoonup}
\newcommand{\ppr}{\rightharpoonup}

\newcommand{\g}{\mathfrak g}

\newcommand{\hg}{\hat{\mathfrak g}}

\newcommand{\h}{V}
\newcommand{\frkg}{\mathfrak g}
\newcommand{\frkh}{V}

\newcommand{\dlam}{_\lambda}

\newcommand{\E}{\mathrm{E}}
\renewcommand{\L}{\mathcal{L}}

\newcommand{\Hom}{\mathrm{Hom}}

\newcommand{\Ker}{\mathrm{Ker}}

\newcommand{\End}{\mathrm{End}}
\newcommand{\Ext}{\mathrm{Ext}}

\newcommand{\ad}{\mathrm{ad}}

\newcommand{\id}{\mathrm{id}}

\begin{document}

\allowdisplaybreaks

\title{Extending structures for 3-Lie algebras}
\author{Tao Zhang}

\date{}
\maketitle



\begin{abstract}
 The cohomology and deformation theory of 3-Lie algebras are revisited.
  The theory of extending structures and unified product for 3-Lie algebras are developed.
  It is proved that the  extending structures of 3-Lie algebras can be classified by using some non-abelian cohomology and deformation map theory.
\end{abstract}

\section{Introduction}
The concept of 3-Lie algebra dates back to Nambu's work  \cite{Nam} to generalize the classical Hamiltonian mechanics.
He derived the following equation
$$\{h_1,h_2,\{x,y,z\}\}=\{\{h_1,h_2,x\},y,z\}+\{x,\{h_1,h_2,y\},z\}+\{x,y,\{h_1,h_2,z\}\},$$
where $h_1,h_2$ are Hamiltonian and $\{\cdot,\cdot,\cdot\}$ is a ternary product.
In \cite{Fil1}, Filippov  introduced the concept of $n$-Lie algebra.
An $n$-Lie algebra is a vector space $\frkg$ with an $n$-ary totally skew-symmetric linear map ($n$-bracket)
from $\bigwedge{}^n\frkg$ to $\frkg$: $(x_1,\cdots, x_n)\mapsto [x_1,\cdots, x_n]$ satisfying the $n$-Jacobi identity or fundamental identity
\begin{eqnarray}\label{eq:Jacobi-n}
[x_1, \cdots, x_{n-1},[y_1,y_2, \cdots , y_n]]=\sum_{i=1}^n [y_1, \cdots,[x_1, \cdots, x_{n-1},y_i],\cdots y_n]
\end{eqnarray}
for all $x_i, y_i \in \frkg$.
When the $n$-ary  linear map is not skew-symmetric, it is called $n$-Leibniz algebra or Leibniz $n$-algebra.
More recently, 3-Lie algebras are applied to the study of gauge symmetry and supersymmetry of multiple coincident M2-branes.


The algebraic theory of $n$-Lie algebras have been studied by many authors, see \cite{DT,Gau,FF,Kas,Ling,Tak2}.
The cohomology theory for $n$-Lie algebras have been studied in \cite{Tak2,CLP,Rot}.

On the other hand, extending structures for Lie algebras, associative algebras and Leibniz algebras are studied by Agore and Militaru in \cite{AM1,AM2,AM3,AM4,AM5}.
It worth also mentioning that the extending structures for left-symmetric algebras, associative and Lie conformal algebras has been studied by Y. Hong and Y. Su in  \cite{Hong1,Hong2,HS}.

This paper is dedicated to study extending structures and unified products for  3-Lie algebras.
We will follow closely to the theory of unified product and extending structures  which are well  developed by A. L. Agore and G. Militaru in \cite{AM1,AM2,AM3}.
Let $\g$ be a 3-Lie algebra and $E$ a vector space containing $\g$ as a subspace.
We will describe and classify all 3-Lie algebras structures on $E$ such that $\g$ is a subalgebra of
$E$. We show that associated to any extending structures of $\g$ by a complement space $V$, there is a  unified product on the direct sum space  $E\cong\g\oplus V$.
We will show how to classify extending structures for  3-Lie algebras by using some non-abelian cohomology and deformation map  theory.

The organization of this paper is as follows.
In Section \ref{sec:2},  we spell out the abelian cohomology theory for 3-Lie algebras.
We  prove that this kind of cohomology theory can be used to  classify abelian extensions of 3-Lie algebras.
We also introduce the notion of Nijenhuis operators for 3-Lie algebras,
which is analogy to the case of ordinary Lie algebras in \cite{D,KM} and of associative algebras in \cite{CGM}.
It is proved that this kind of operator gives trivial deformation.
Note that some results of this section has been generalized to the cases of Lie triple systems and 3-Lie colour algebras  in \cite{Zhang1, Zhang2}.
In Section \ref{sec:extending},  we study extending structures for 3-Lie algebras.
We will describe and classify all 3-Lie algebras structures on $E$ such that $\g$ is a subalgebra of
$E$. We show that associated to any extending structures of $\g$ by a complement space $V$, there is a  unified product on the direct sum space  $E\cong\g\oplus V$.
In Section \ref{sec:specialcases}, we give two special cases of the
unified product, namely the crossed product and bicrossed product   of 3-Lie algebras
which are related to the study of the non-abelian extension problem  and classifying complements problem respectively.
It is proved that  they can be classified by using some non-abelian cohomology and deformation map  theory.

Throughout this paper, all vector spaces are assumed to be over
an algebraically closed field of characteristic not equal to 2 and 3.
Let $V$ and $W$ be two vector spaces. The space of linear maps from $V$ to $W$ is denoted by $\Hom(V,W)$.

\section{3-Lie algebras and cohomology}\label{sec:2}

In this section, we give a review on the abelian cohomology and deformation theory of 3-Lie algebras, the main references are \cite{Bai,CLP,LSZB,LZ,Tak2,Zhang2}.

A 3-Lie algebra consists of a vector space $\frkg$ together with a bracket $[\cdot,\cdot,\cdot]: \bigwedge^3 \g\to \g$  such that the fundamental identity
\begin{eqnarray}\label{eq:Jacobi-3}
[x_1, x_2, [y_1, y_2, y_3]] = [[ x_1, x_2, y_1], y_2, y_3] + [y_1, [ x_1, x_2, y_2], y_3] + [y_1, y_2, [ x_1, x_2, y_3]]
\end{eqnarray}
holds for all $x_i, y_i \in \g$.

Since the bracket $[\cdot,\cdot,\cdot]$ is skew-symmetric, the above fundamental identity \eqref{eq:Jacobi-3}
can be rewritten as
\begin{eqnarray}\label{eq:Jacobi-cp}
[x_1, x_2, [y_1, y_2, y_3]] &=& [[ x_1, x_2, y_1], y_2, y_3] + [[ x_1, x_2, y_2], y_3, y_1] + [[ x_1, x_2, y_3],y_1, y_2]\\
&:=& [[ x_1, x_2, y_1], y_2, y_3] +c.p.
\end{eqnarray}
where c.p. means cyclic permutations with respect to elements $y_1, y_2$ and $y_3$.

Denote by $x=(x_1, x_2)$ and $\ad(x) y_i=[x_1, x_2, y_i]$, then the above equality can be rewritten in the form
\begin{eqnarray}\label{eq:Jacobi-3'}
\ad(x)[y_1, y_2, y_3]=[\ad(x)y_1, y_2, y_3]+[y_1, \ad(x)y_2, y_3]+[y_1, y_2, \ad(x) y_3].
\end{eqnarray}

Denoted by $\L:=\bigwedge^2 \frkg$, which is called fundamental set.
The elements $x=(x_1,x_2)\in \bigwedge^2 \frkg$ are called fundamental object.
Define an operation on fundamental object by
\begin{eqnarray}\label{eq:fundamental}
x\circ y=([x_1,x_2,y_1],y_2)+(y_1,[x_1,x_2,y_2]).
\end{eqnarray}
In \cite{DT}, the authors  proved that $\L$ is a Leibniz algebra  satisfying the following Leibniz identity
$$x\circ (y\circ z)= (x\circ y)\circ z+y\circ (x\circ z),$$
and
$$\ad(x)\ad(y) w-\ad(y)\ad(x) w=\ad(x\circ y) w,$$ 
for all $x,y,z\in\L, w\in \frkg$. Thus $\ad: \L \to \End(\frkg)$ is a homomorphism of Leibniz algebras.

Recall that for a Leibniz algebra $\L$, a representation is a vector space $V$ together with two bilinear maps
$$[\cdot,\cdot]_L:\L\times V\to V \,\,\, \text{and}\,\,\,  [\cdot,\cdot]_R: V\times \L\to V$$
satisfying the following three axioms
\begin{itemize}
\item[$\bullet$] {\rm(LLM)}\quad  $[x\circ y, m]_L=[x, [y, m]_L]_L-[y,[x, m]_L]_L$,
\item[$\bullet$] {\rm(LML)}\quad  $[m, x\circ y]_R=[[m, x]_R, y]_R+[x, [m, y]_R]_L$,
\item[$\bullet$] {\rm(MLL)}\quad  $[m, x\circ y]_R=[x,[m, y]_R]_L-[[x,m]_L, y]_R.$
\end{itemize}
By (LML) and (MLL) we also have
\begin{itemize}
\item[$\bullet$] {\rm(MMM)}\quad  $[[m, x]_R, y]_R+[[x,m]_L, y]_R=0.$
\end{itemize}
In fact, assume (LLM), one of (LML),(MLL),(MMM) can be derived from the other two.

Given a 3-Lie algebra $\frkg$ and a vector spaces $\frkh$, we define the maps
$$[\cdot,\cdot]_L:\L \otimes \Hom(\frkg,\frkh)\to \Hom(\frkg,\frkh)\quad \mbox{and} \quad [\cdot,\cdot]_R: \Hom(\frkg,\frkh)\otimes \L \to \Hom(\frkg,\frkh)$$
by
\begin{eqnarray}
\label{eq:leibniz01}{[(x_1,x_2),\phi]_L}(x_3)&=&\rho(x_1,x_2)\phi(x_3)-\phi([x_1,x_2,x_3]),\\
\nonumber{[\phi,(x_1,x_2)]_R}(x_3)&=&\phi([x_1,x_2,x_3])-\rho(x_1, x_2)\phi(x_3)-\rho(x_2,x_3)\phi(x_1)\\
\label{eq:leibniz02}&&-\rho(x_3,x_1)\phi(x_2),
\end{eqnarray}
for all $\phi\in \Hom(\frkg,\frkh), x_i\in \frkg$, where $\rho$ is a map from $\L=\bigwedge{}^2\frkg$ to $\End(V)$.

\begin{pro}\label{pro:rep}(\cite{CLP, Zhang2})
Let $\frkg$ be a 3-Lie algebra. Then $\Hom(\frkg,\frkh)$ equipped with the above two maps
$[\cdot,\cdot]_L$ and $[\cdot,\cdot]_R$ is a representation of the Leibniz algebra $\L$
if and only if the following two conditions are satisfied for all $x_i, y_i \in \frkg$,
\begin{itemize}
\item[$\bullet$]{\rm(R1)}\quad $[\rho(x_1,x_2),\rho(y_1,y_2)]=\rho((x_1,x_2)\circ (y_1,y_2))$,
\item[$\bullet$]{\rm(R2)}\quad $\rho(x_1,[y_1, y_2, y_3])=\rho(y_1, y_2)\rho(x_1,y_3)+\rho(y_2, y_3)\rho(x_1,y_1) + \rho(y_3, y_1)\rho(x_1,y_2)$.
\end{itemize}
\end{pro}

For the proof of the above Proposition \ref{pro:rep}, see \cite{CLP, Zhang2}.

\begin{definition} Let $\frkg$ be a 3-Lie algebra and $V$ be a vector space. Then $(V, \rho)$ is called a representation
of $\frkg$  if and only if the conditions (R1) and (R2) in the above Proposition \ref{pro:rep} are satisfied.
\end{definition}

If we denote by $(x_1,x_2)\trr u:= \rho(x_1,x_2)(u)$, then  the conditions (R1) and (R2) can be rewritten as
\begin{itemize}
\item[$\bullet$]{\rm(R1')}\quad $[(x_1,x_2),(y_1,y_2)]\trr u=\big((x_1,x_2)\circ (y_1,y_2)\big)\trr u$,
\item[$\bullet$]{\rm(R2')}\quad $(x_1,[y_1, y_2, y_3])\trr u=(y_1, y_2)\trr \big((x_1,y_3)\trr u\big)+(y_2, y_3)\trr \big((x_1,y_1)\trr u\big) + (y_3, y_1)\trr \big((x_1,y_2)\trr u\big)$.
\end{itemize}

For example, given a 3-Lie algebra $\frkg$, there is a natural \emph{adjoint representation} on itself. The corresponding representation $\ad(x_1,x_2)$ is given by
\begin{eqnarray*}
\ad(x_1,x_2)(x_3):=[x_1,x_2,x_3].
\end{eqnarray*}

Now we define the cochain complex for a 3-Lie algebra $\frkg$ with coefficients in $\frkh$ by
$$C^n(\frkg,\frkh):=\Hom\left(\bigwedge{}^{2n+1}\frkg,\frkh\right)\subseteq\Hom\left(\left(\bigwedge{}^{2n}\frkg\right)\otimes\frkg,\frkh\right)
\cong\Hom\left(\L{}^{n},\Hom(\frkg,\frkh)\right)$$
and
$$d_{n-1}:C^{n-1}(\frkg,\frkh)\to C^{n}(\frkg,\frkh)$$
where
\begin{eqnarray*}
&&d_{n-1}\omega(x^1,x^2,\cdots,x^n,w)\\
&=&d_{n-1}\omega(x^1,x^2,\cdots,x^n)(w)\\
&=&\sum_{k=1}^{n-1}(-1)^{k+1}[x^k,\omega(x^1,\cdots,\hat{x^k},\cdots,x^n)]_L(w)+(-1)^n[\omega(x^1,\cdots,x^{n-1}),x^n]_R(w)\\
&&+\sum_{1\leq k<j\leq n}(-1)^{k}\omega(x^1,\cdots,x^{j-1},x^k\circ x^j,x^{j+1},\cdots,x^n)(w),
\end{eqnarray*}
for all $x^i\in \L=\bigwedge{}^2\frkg,\ w\in \frkg$.
In other words, we define the cohomology of a 3-Lie algebra $\frkg$ with coefficients in $\frkh$ to be
the cohomology of Leibniz algebra $\L$ with coefficients in $\Hom(\frkg,\frkh)$.
For more details of cohomology of  Leibniz algebras, see \cite{Loday}.

\begin{thm}(\cite{Zhang2})
Let $\frkg$ be a 3-Lie algebra and $(V,\rho)$ be a representation of $\frkg$.
Then there exists a cochain complex $\left\{C(\frkg,\frkh)=\bigoplus_{n\geq 0}C^n(\frkg,\frkh),d\right\}$,
where the coboundary operator is given by
\begin{eqnarray}\label{eq:coboudary}
\notag&&d_{n-1}\omega (x_1, x_2,\cdots, x_{2n+1})\\
\notag&=&(-1)^{n+1}\rho(x_{2n+1}, x_{2n-1})\omega (x_1, x_2,\cdots, x_{2n-2}, x_{2n})\\
\notag&&+(-1)^{n+1}\rho(x_{2n}, x_{2n+1})\omega (x_1, x_2,\cdots, x_{2n-1})\\
\notag&&+\sum_{k=1}^n (-1)^{k+1}\rho(x_{2k-1}, x_{2k})\omega (x_1, x_2,\cdots,\widehat{x_{2k-1}}, \widehat{x_{2k}},\cdots, x_{2n+1})\\
&&+\sum_{k=1}^n\sum^{2n+1}_{j=2k+1}(-1)^{k}
\omega(x_1, x_2,\cdots,\widehat{x_{2k-1}}, \widehat{x_{2k}},\cdots, [x_{2k-1}, x_{2k}, x_j],\cdots, x_{2n+1})
\end{eqnarray}
such that $d\circ d=0$.
\end{thm}

\pf
Put $x^k=(x_{2k-1}, x_{2k})$, $w=x_{2n+1}$ and $[\cdot,\cdot]_L$, $[\cdot,\cdot]_R$ as in \eqref{eq:leibniz01} and \eqref{eq:leibniz02},
then we get a coboundary operator $d_{n-1}: C^{n-1}(\g, V)\to C^{n}(\g, V)$ as the following:
\begin{eqnarray*}
&&d_{n-1} \omega (x_1, x_2,\cdots, x_{2n+1}):=d_{n-1}\omega(x^1,x^2,\cdots,x^n)(w)\\
&=&\sum_{k=1}^{n-1} (-1)^{k+1}[x_{2k-1}, x_{2k},\omega (x_1, x_2,\cdots,\widehat{x_{2k-1}}, \widehat{x_{2k}},\cdots, x_{2n})]_L(x_{2n+1})\\
&&+(-1)^n[\omega(x_1,x_2,\cdots,x_{2n-2}),x_{2n-1}, x_{2n}]_R(x_{2n+1})\\
&&+\sum_{k=1}^n\sum^{2n}_{j=2k+1}(-1)^{k}\omega(x_1,\cdots,\widehat{x_{2k-1}}, \widehat{x_{2k}},\cdots,[x_{2k-1}, x_{2k}, x_j],\cdots, x_{2n})(x_{2n+1})\\
&=&\sum_{k=1}^{n-1} (-1)^{k+1}\{\rho(x_{2k-1}, x_{2k})\omega (x_1, x_2,\cdots,\widehat{x_{2k-1}}, \widehat{x_{2k}},\cdots, x_{2n+1})\\
&&-\omega(x_1, x_2,\cdots,\widehat{x_{2k-1}}, \widehat{x_{2k}},\cdots, [x_{2k-1}, x_{2k}, x_{2n+1}])\}\\[1em]
&&+(-1)^n \{\omega(x_1, x_2,\cdots,[x_{2n-1}, x_{2n}, x_{2n+1}])-\rho(x_{2n-1}, x_{2n})\omega (x_1, x_2,\cdots, x_{2n+1})\\
&&-\rho(x_{2n+1}, x_{2n-1})\omega (x_1, x_2,\cdots, x_{2n-2}, x_{2n})-\rho(x_{2n},x_{2n+1})\omega (x_1, x_2,\cdots, x_{2n-1})\}\\[1em]
&&+\sum_{k=1}^n\sum^{2n}_{j=2k+1}(-1)^{k}
\omega(x_1, x_2,\cdots,\widehat{x_{2k-1}}, \widehat{x_{2k}},\cdots, [x_{2k-1}, x_{2k}, x_j],\cdots,x_{2n+1})\\
&=&(-1)^{n+1}\rho(x_{2n+1}, x_{2n-1})\omega (x_1, x_2,\cdots, x_{2n-2}, x_{2n})\\
&&+(-1)^{n+1}\rho(x_{2n}, x_{2n+1})\omega (x_1, x_2,\cdots, x_{2n-1})\\
&&+\sum_{k=1}^n (-1)^{k+1}\rho(x_{2k-1}, x_{2k})\omega (x_1, x_2,\cdots,\widehat{x_{2k-1}}, \widehat{x_{2k}},\cdots, x_{2n+1})\\
&&+\sum_{k=1}^n\sum^{2n+1}_{j=2k+1}(-1)^{k}
\omega(x_1, x_2,\cdots,\widehat{x_{2k-1}}, \widehat{x_{2k}},\cdots, [x_{2k-1}, x_{2k}, x_j],\cdots, x_{2n+1}).
\end{eqnarray*}

\qed

\begin{definition}
The quotient space
${\mathcal{H}}^{\bullet}(\frkg, V)=Z^{\bullet}(\frkg, V)/B^{\bullet}(\frkg, V)$,
where $Z^{\bullet}(\frkg, V)=\{\omega\in C^{n}(\g, V)|d\omega=0\}$ is the space of cocycles
and $B^{\bullet}(\frkg, V)=\{\omega=d\nu|\nu\in C^{n-1}(\g, V)\}$ is the space of coboundaries,
is called the cohomology group of a 3-Lie algebra $\frkg$ with coefficients in $\frkh$.
\end{definition}

According to the above definition, a 0-cochain is a map $\nu\in\Hom(\frkg,\frkh)$,
a 1-cochain is a map
$\omega\in \Hom\left(\bigwedge{}^3\frkg,\frkh\right)$,
and the coboundary operator is give by
\begin{eqnarray}
\label{eq:cobound01}d_{0}\nu(x^1,w)&=&d_{0}\nu(x^1)(w)=-[\nu,x^1]_R(w),\\
\label{eq:cobound02}d_{1}\omega(x^1,x^2,w)&=&[x^1,\omega(x^2)]_L(w)+[\omega(x^1),x^2]_R(w)-\omega(x^1\circ x^2)(w).
\end{eqnarray}

Put $x^1=(x_1,x_2)\in \L,\ w=x_3\in\frkg$ in the equality \eqref{eq:cobound01}, then by \eqref{eq:leibniz02} we have
\begin{eqnarray}
\nonumber&&d_0\nu(x_1, x_2,x_3)\\
&=&\rho(x_1, x_2)\nu(x_3)+\rho(x_1,x_3)\nu(x_2)+\rho(x_2,x_3)\nu(x_1)-\nu([x_1, x_2, x_3]).
\end{eqnarray}

\begin{definition}
Let $\frkg$ be a 3-Lie algebra and $(V, \rho)$ be a representation of $\frkg$. Then a map $\nu\in\Hom(\frkg,\frkh)$
is called 0-cocycle if and only if $\forall\, x_1, x_2, x_3\in \frkg$,
\begin{eqnarray}\label{eq:0coc}
\rho(x_1, x_2)\nu(x_3)+\rho(x_1,x_3)\nu(x_2)+\rho(x_2,x_3)\nu(x_1)-\nu([x_1, x_2, x_3])=0,
\end{eqnarray}
and a map $\omega\in \Hom\left(\bigwedge{}^3\frkg,\frkh\right)$ is called a 1-coboudary if there exists a map $\nu\in\Hom(\frkg,\frkh)$ such that $\omega=d_0\nu$.
\end{definition}

Put $x^1=(x_1,x_2)\in \L, x^2=(y_1,y_2)\in \L$, $w=y_3\in\frkg$ in the equality \eqref{eq:cobound02}, then we have
\begin{eqnarray*}
d_{1}\omega(x_1,x_2,y_1,y_2,y_3)&=&[x_1,x_2,\omega(y_1,y_2)]_L(y_3)+[\omega(x_1,x_2),y_1,y_2]_R(y_3)\\
&&-\omega((x_1,x_2)\circ (y_1,y_2))(y_3),
\end{eqnarray*}
where
\begin{eqnarray*}
[x_1,x_2,\omega(y_1,y_2)]_L(y_3)&=&\rho(x_1,x_2)\omega(y_1,y_2)(y_3)-\omega(y_1,y_2)([x_1,x_2,y_3])\\
&=&\rho(x_1,x_2)\omega(y_1,y_2,y_3)-\omega(y_1,y_2,[x_1,x_2,y_3]),
\end{eqnarray*}
\begin{eqnarray*}
[\omega(x_1,x_2),y_1,y_2]_R(y_3)&=&\omega(x_1,x_2)([y_1,y_2,y_3])-\rho(y_1,y_2)\omega(x_1,x_2)(y_3)\\
&&-\rho(y_1,y_3)\omega(x_1,x_2)(y_2)-\rho(y_2,y_3)\omega(x_1,x_2)(y_1)\\
&=&\omega(x_1,x_2,[y_1,y_2,y_3])-\rho(y_1,y_2)\omega(x_1,x_2,y_3)\\
&&-\rho(y_1,y_3)\omega(x_1,x_2,y_2)-\rho(y_2,y_3)\omega(x_1,x_2,y_1),
\end{eqnarray*}
\begin{eqnarray*}
\omega((x_1,x_2)\circ (y_1,y_2))(y_3)&=&\omega(([x_1,x_2,y_1],y_2)+(y_1,[x_1,x_2,y_2]))(y_3)\\
&=&\omega(([x_1,x_2,y_1],y_2,y_3)+\omega(y_1,[x_1,x_2,y_2],y_3).
\end{eqnarray*}
\begin{definition}\label{def:1coc}
Let $\frkg$ be a 3-Lie algebra and $(V, \rho)$ be a representation of $\frkg$. Then a map
$\omega\in \Hom\left(\bigwedge{}^3\frkg,\frkh\right)$ is called 1-cocycle if and only if $\forall\, x_1, x_2,y_1, y_2, y_3\in \frkg$,
\begin{eqnarray}\label{eq:1coc}
\nonumber&& \omega(x_1, x_2,[y_1, y_2, y_3])+\rho(x_1, x_2)\omega(y_1, y_2, y_3)\\
\nonumber&=&\omega([x_1, x_2, y_1], y_2, y_3) + \omega([x_1, x_2,y_2], y_3, y_1) + \omega(y_1, y_2, [x_1, x_2,y_3])\\
&&+\rho( y_2, y_3)\omega(x_1,x_2,y_1) + \rho(y_3, y_1)\omega(x_1,x_2,y_2) + \rho(y_1, y_2)\omega(x_1,x_2,y_3).
\end{eqnarray}
\end{definition}


Now we give a brief description of  the theory of abelian extensions of 3-Lie algebras.
We are going to show that associated to any abelian extension, there is a representation and a 1-cocycle.
Furthermore, abelian extensions can be classified by the first cohomology group.

\begin{definition}
 Let $(\frkg, [\cdot,\cdot,\cdot]_\g)$, $(\frkh, [\cdot,\cdot,\cdot]_\frkh)$,
 $(\hat{\frkg}, [\cdot,\cdot,\cdot]_{\hat{\g}})$ be 3-Lie algebras and
$i:\frkh\to\hat{\frkg},~~p:\hat{\frkg}\to\frkg$
be homomorphisms. The following sequence of 3-Lie algebras is a
short exact sequence if $\mathrm{Im}(i)=\mathrm{Ker}(p)$,
$\mathrm{Ker}(i)=0$ and $\mathrm{Im}(p)=\g$,
\begin{equation}\label{diagram:exact}
 \xymatrix{
   0  \ar[r]^{} & \frkh \ar[r]^{i} & \hat{\frkg} \ar[r]^{p} & \frkg  \ar[r]^{} & 0. \\
  }
\end{equation}
In this case, we call $\hat{\frkg}$  an extension of $\frkg$ by
$\frkh$, and denote it by $\E_{\hat{\g}}$.
It is called an abelian extension if $\frkh$ is abelian ideal of $\hat{\frkg}$, i.e. $[u,v,\cdot]_{\hat{\g}}=[u,\cdot,v]_{\hat{\g}}=[\cdot,u,v]_{\hat{\g}}=0$, $\forall u,v\in \frkh$.
\end{definition}

A section $\sigma:\frkg\to\hat{\frkg}$ of $p:\hat{\frkg}\to\frkg$
consists of linear maps
$\sigma:\frkg\to\hat{\g}$
 such that  $p\circ\sigma=\id_{\frkg}$.

\begin{definition}
 Two extensions of 3-Lie algebras
 $\E_{\hat{\g}}:0\to\frkh\stackrel{i}{\to}\hat{\g}\stackrel{p}{\to}\g\to0$
 and $\E_{\tilde{\g}}:0\to\frkh\stackrel{j}{\to}\tilde{\g}\stackrel{q}{\to}\g\to0$ are equivalent,
 if there exists a 3-Lie algebra homomorphism $F:\hat{\g}\to\tilde{\g}$  such that the following diagram commutes
\begin{equation}\label{diagram:equivalent}
\xymatrix{
   0  \ar[r]^{} & \frkh \ar[d]_{\id_V} \ar[r]^{i} & \hat{\frkg} \ar[d]_{F} \ar[r]^{p} & \frkg \ar[d]_{\id_{\g}} \ar[r]^{} & 0 \\
   0 \ar[r]^{} & \frkh \ar[r]^{j} & \tilde{\frkg} \ar[r]^{q} & \frkg \ar[r]^{} & 0
   }
\end{equation}
The set of equivalent classes of extensions of $\g$ by $\h$ is denoted by $\Ext(\frkg,\frkh)$.
\end{definition}

Let $\hat{\frkg}$ be an abelian extension of $\frkg$ by
$\frkh$, and $\sigma:\g\to\hat{\g}$ be a section.  Denote by
$$\si(x)=\si(x_1,x_2)=(\si x_1,\si x_2),$$
and define $\rho:\wedge^2\frkg\to\End(\frkh)$ by
\begin{equation}\label{eq:rep}
\rho(x)(u)=\rho(x_1,x_2)(u)\triangleq[\sigma(x_1),\sigma(x_2),u]_{\hat{\g}}=\ad(\si(x))u,
\end{equation}
for all $x=(x_1,x_2)\in\bigwedge^2\frkg$, $u\in\frkh$.

\begin{lem}\label{pro:2-modules}
With the above notations, $\rho$ is a representation of $\g$ on $\h$ and does not depend on the choice of the section $\sigma$.
Moreover,  equivalent abelian extensions give the same representation of $\g$ on $\h$.
\end{lem}

\pf We give a sketch of the proof.
First, we show that $\rho$ is independent of
the choice of $\sigma$. In fact, if we choose another section $\sigma':\g\to\hg$, then
$$p(\sigma(x_i)-\sigma'(x_i))=x_i-x_i=0
\Longrightarrow\sigma(x_i)-\sigma'(x_i)\in \h\Longrightarrow\sigma'(x_i)=\sigma'(x_i)+u$$
for some $u\in\h$.

Since we  have $[\cdot,u,v]_{\hat{\g}}=0$
for all $u,v\in\frkh$, which implies that
\begin{eqnarray*}
[\sigma'(x_1),\sigma'(x_2),w]_{\hat{\g}}&=&[\sigma(x_1)+u,\sigma(x_2)+v,w]_{\hat{\g}}\\
&=&[\sigma(x_1),\sigma(x_2)+v,w]_{\hat{\g}}+[u,\sigma(x_2)+v,w]_{\hat{\g}}\\
&=&[\sigma(x_1),\sigma(x_2),w]_{\hat{\g}}+[\sigma(x_1),v,w]_{\hat{\g}}\\
&=&[\sigma(x_1),\sigma(x_2),w]_{\hat{\g}},
\end{eqnarray*}
thus $\rho$ is independent on the choice of $\sigma$.

Second, one can show that $\rho$ is a representation of $\g$ on $\h$ since $\h$ is an abelian ideal of $\g$.
%
%

At last, suppose that $\E_{\hat{\g}}$ and $\E_{\tilde{\g}}$ are equivalent abelian extensions, and $F:\hat{\g}\to\tilde{\g}$ is the 3-Lie algebra homomorphism satisfying $F\circ i=j$, $q\circ F=p$.
Choose linear sections $\sigma$ and $\sigma'$ of $p$ and $q$, we get $qF\sigma(x_i)=p\sigma(x_i)=x_i=q\sigma'(x_i)$,
then $F\sigma(x_i)-\sigma'(x_i)\in \Ker (q)\cong\h$. Thus, we have
$$
[\sigma(x_1),\sigma(x_2),u]_{\hat{\g}}=[F\sigma(x_1),F\sigma(x_2),u]_{\tilde{\g}}=[\sigma'(x_1),\sigma'(x_2),u]_{\tilde{\g}}.
$$
Therefore, equivalent abelian extensions give the same $\rho$. The proof is finished.
\qed\vspace{3mm}

Let $\sigma:\frkg\to\hat{\frkg}$  be a
section of the abelian extension. Define the following map:
\begin{equation}\label{eq:coc}
\omega(x_1,x_2,x_3)\triangleq[\sigma(x_1),\sigma(x_2),\sigma(x_3)]_{\hat{\frkg}}-\sigma([x_1,x_2,x_3]_\frkg),
\end{equation}
for all $x_1,x_2,x_3\in\frkg$. By direct computations,  we have

\begin{lem}\label{thm:2-cocylce}
Let $0\to\frkh{\to}\hat{\g}{\to}\g\to 0$ be an abelian extension of $\g$ by $\h$.
Then $\omega$ defined by \eqref{eq:coc} is a 1-cocycle of $\g$ with coefficients in $\frkh$,
where the representation $\rho$ is given by \eqref{eq:rep}.
\end{lem}

Now we can transfer the 3-Lie algebra structure on $\hat{\g}$ to the 3-Lie algebra structure on $\frkg\oplus\frkh$ using the 1-cocycle given above. The following Lemma \ref{lem:prod} and Lemma \ref{lem:2-cocylce}  can be proved directly,  see \cite{Bai,LZ} for more details.

\begin{lem}\label{lem:prod}
Let $\frkg$ be a 3-Lie algebra, $(V, \rho)$ be an $\frkg$-module and $\omega: \bigwedge^3 \g\to V$ is a 1-cocycle.
Then $\g\oplus V$ is a 3-Lie algebra under the following bracket:
\begin{eqnarray*}
&&[x_1 + u_1, x_2 + u_2, x_3 + u_3]_{\omega}\\
&=&[x_1, x_2, x_3] + \omega(x_1, x_2, x_3)+ \rho (x_1, x_2)(u_3) + \rho( x_2, x_3)(u_1)+ \rho( x_3, x_1)(u_2),
\end{eqnarray*}
where $x_1, x_2, x_3 \in \frkg$ and $u_1, u_2, u_3 \in V$. This 3-Lie algebra is denoted by $\g\oplus_{\omega}\h$.
\end{lem}

\begin{lem}\label{lem:2-cocylce}
 Two abelian extensions of 3-Lie algebras $0\to\frkh{\to}\g\oplus_\omega\h{\to}\g\to 0$
 and  $0\to\frkh{\to}\g\oplus_{\omega'}\h{\to}\g\to 0$ are equivalent if and only if $\omega$ and $\omega'$ are in the same cohomology class.
\end{lem}

From the above Lemmas \ref{pro:2-modules}--\ref{lem:2-cocylce},  we obtain
\begin{thm}
Let $\frkg$ be a 3-Lie algebra and $(V,\rho)$  a representation of $\g$.
Then there is a one-to-one correspondence between equivalence classes  $\Ext(\frkg,\frkh)$ of abelian extensions of the 3-Lie algebra $\frkg$ by $\h$ and the first cohomology group ${\mathcal{H}}^1(\frkg,\h)$.
\end{thm}


Next,  we study infinitesimal deformations and Nijenhuis operators for 3-Lie algebras.
This kind of operator gives trivial deformation.

Let $\frkg$ be a 3-Lie algebra and $\omega:\bigwedge{}^3\g\to\g$ be a linear map. Consider a $\lambda$-parametrized family of linear maps:
\begin{eqnarray*}
[x_1, x_2, x_3]_\lam&\triangleq& [x_1, x_2, x_3]+ \lambda\omega(x_1, x_2, x_3).
 \end{eqnarray*}

If $[\cdot,\cdot,\cdot]_\lam$ endow $\frkg$ with a 3-Lie algebra structure, then we say that $\omega$ generates a
first order deformation of the 3-Lie algebra $\frkg$.


\begin{thm}\label{thm:deformation}
$\omega$ generates a first order deformation of the 3-Lie algebra $\frkg$ is equivalent to
(i) $\omega$ itself defines a 3-Lie algebra structure on $\g$ and
(ii) $\omega$ is a 1-cocycle of $\frkg$ with the coefficients in the adjoint representation.
\end{thm}

\pf
For the equality
\begin{eqnarray*}
&&[x_1, x_2, [y_1, y_2, y_3]\dlam]\dlam \\
&=&[[ x_1, x_2, y_1]\dlam , y_2, y_3]\dlam  + [y_1, [ x_1, x_2, y_2]\dlam , y_3]\dlam  + [y_1, y_2, [ x_1, x_2, y_3]\dlam]\dlam,
\end{eqnarray*}
the left hand side is equal to
\begin{eqnarray*}
&&[x_1, x_2, [y_1, y_2, y_3]+\lam\omega(y_1, y_2, y_3)]\dlam\\
&=&[x_1, x_2, [y_1, y_2, y_3]]+\lam\omega(x_1, x_2, [y_1, y_2, y_3])\\
&&+[x_1, x_2,\lam\omega(y_1, y_2, y_3)]+\lam\omega(x_1, x_2,\lam\omega(y_1, y_2, y_3))\\
&=&[x_1, x_2, [y_1, y_2, y_3]]+\lam\{\omega(x_1, x_2, [y_1, y_2, y_3])+[x_1, x_2,\omega(y_1, y_2, y_3)]\}\\
&&+\lam^2\omega(x_1, x_2,\omega(y_1, y_2, y_3)),
\end{eqnarray*}
and the right hand side is equal to
\begin{eqnarray*}
&&[[x_1, x_2, y_1]+\lam\omega(x_1, x_2, y_1), y_2, y_3]\dlam  + [y_1, [ x_1, x_2, y_2]+\lam\omega(x_1, x_2, y_2), y_3]\dlam \\
&&+ [y_1, y_2, [ x_1, x_2, y_3]+\lam\omega(x_1, x_2, y_3)]\dlam\\
&=&[[ x_1, x_2, y_1] , y_2, y_3]  + [y_1, [ x_1, x_2, y_2] , y_3] + [y_1, y_2, [ x_1, x_2, y_3]]\\
&&+\lam\{\omega([x_1, x_2, y_1],y_2, y_3) +[\omega(x_1, x_2, y_1), y_2, y_3]\\
&&\qquad+\omega(y_1,[ x_1, x_2, y_2], y_3)+[y_1,\omega(x_1, x_2, y_2), y_3]\\
&&\qquad+\omega(y_1, y_2, [ x_1, x_2, y_3])+[y_1,y_2,\omega(x_1, x_2, y_3)] \}\\
&&+\lam^2\{\omega(\omega( x_1, x_2, y_1), y_2, y_3) + \omega(y_1, \omega( x_1, x_2, y_2), y_3) + \omega(y_1, y_2, \omega( x_1, x_2, y_3))\}.
\end{eqnarray*}
Thus we have
\begin{eqnarray}
\nonumber &&\omega(x_1, x_2, [y_1, y_2, y_3])+[x_1, x_2,\omega(y_1, y_2, y_3)]\\
\nonumber &=&\omega([x_1, x_2, y_1],y_2, y_3)+\omega(y_1,[ x_1, x_2, y_2], y_3) +\omega(y_1, y_2, [ x_1, x_2, y_3])\\
\label{eq:2-coc01} &&+[\omega(x_1, x_2, y_1), y_2, y_3]+[y_1,\omega(x_1, x_2, y_2), y_3]+[y_1,y_2,\omega(x_1, x_2, y_3)],\\
\nonumber &&\omega(x_1, x_2,\omega(y_1, y_2, y_3))\\
\label{eq:2-coc02}&=&\omega(\omega( x_1, x_2, y_1), y_2, y_3) + \omega(y_1, \omega( x_1, x_2, y_2), y_3)+ \omega(y_1, y_2, \omega( x_1, x_2, y_3)).
\end{eqnarray}
Therefore $\omega$ defines a 3-Lie algebra structure on $\g$ and $\omega$ is a 1-cocycle of $\frkg$ with the coefficients in the adjoint representation.
\qed

\emptycomment{
Two deformations $\omega$ and $\omega'$ are called to be equivalent if there exists an isomorphism
$F=\id+\lambda N:\g_\omega\to\g_{\omega'}$ such that
\begin{eqnarray}
T_\lam [x_1,x_2,x_3]_\omega=[T_\lam x_1,T_\lam x_2,T_\lam x_3]_{\omega'}.
\end{eqnarray}
\begin{eqnarray*}
T_\lam [x_1,x_2,x_3]_\omega&=&[x_1,x_2,x_3]+\lam\omega(x_1,x_2,x_3)+\lam N([x_1,x_2,x_3]+\lam\omega(x_1,x_2,x_3))\\
&=&[x_1,x_2,x_3]+\lam(\omega(x_1,x_2,x_3)+ N[x_1,x_2,x_3])+\lam^2N\omega(x_1,x_2,x_3),
\end{eqnarray*}
\begin{eqnarray*}
[T_\lam x_1,T_\lam x_2,T_\lam x_3]_{\omega'}&=&[x+\lam Nx_1,x_2+\lam Nx_2,x_3+\lam Nx_3]+\lambda{\omega'}(x_1+\lam Nx_1,x_2+\lam Nx_2,x_3+\lam Nx_3)\\
&=&[x_1,x_2,x_3]+\lam([Nx_1,x_2,x_3]+[x_1,Nx_2,x_3]+[x_1,x_2,Nx_3]+\omega'(x_1,x_2,x_3))\\
&&+\lam^2([Nx_1,Nx_2,x_3]+[Nx_1,x_2,Nx_3]+[x_1,Nx_2,Nx_3])\\
&&+\lam^3[Nx_1,Nx_2,Nx_3].
\end{eqnarray*}
}


Recently,  Nijenhuis operator for $n$-Lie algebras is introduced in \cite{LSZB}.
The following results can be derived as special cases when $n=3$ from that paper.

In this case we consider the second order deformation $\frkg_\lam$ by
\begin{eqnarray*}
[x_1, x_2, x_3]_\lam&\triangleq& [x_1, x_2, x_3]+ \lambda\omega_1(x_1, x_2, x_3)+ \lambda^2\omega_2(x_1, x_2, x_3).
 \end{eqnarray*}
A second order deformation is said to be \emph{trivial} if there exists a linear map $N:\frkg\to \frkg$
such that for $T_\lam = \id + \lambda N$: $\frkg_\lam \to \frkg$ there holds
\begin{eqnarray}
T_\lam ([x_1,x_2,x_3]\dlam)=[T_\lam x_1,T_\lam x_2,T_\lam x_3].
\end{eqnarray}
By definition the left hand side equals to
\begin{eqnarray*}
T_\lam ([x_1,x_2,x_3]\dlam)&=&[x_1,x_2,x_3]+\lam(\omega_1(x_1,x_2,x_3)+ N[x_1,x_2,x_3])\\
&&+\lam^2(\omega_2(x_1,x_2,x_3)+N\omega_1(x_1,x_2,x_3))+\lam^3N\omega_2(x_1,x_2,x_3),
\end{eqnarray*}
and the right hand side equals to
\begin{eqnarray*}
[T_\lam x_1,T_\lam x_2,T_\lam x_3]&=&[x_1+\lam Nx_1,x_2+\lam Nx_2,x_3+\lam Nx_3]\\
&=&[x_1,x_2,x_3]+\lam([Nx_1,x_2,x_3]+[x_1,Nx_2,x_3]+[x_1,x_2,Nx_3])\\
&&+\lam^2([Nx_1,Nx_2,x_3]+[Nx_1,x_2,Nx_3]+[x_1,Nx_2,Nx_3])\\
&&+\lam^3[Nx_1,Nx_2,Nx_3].
\end{eqnarray*}
Thus we have
\begin{eqnarray}
\label{eq:Nijenhuis1}&&\omega_1(x_1,x_2,x_3)+N[x_1,x_2,x_3]=[Nx_1,x_2,x_3]+[x_1,Nx_2,x_3]+[x_1,x_2,Nx_3],\\
\notag&&\omega_2(x_1,x_2,x_3)+N\omega_1(x_1,x_2,x_3)\\
\label{eq:Nijenhuis2}&&\qquad\qquad\qquad\quad=[Nx_1,Nx_2,x_3]+[Nx_1,x_2,Nx_3]+[x_1,Nx_2,Nx_3],\\
\label{eq:Nijenhuis3}&&N\omega_2(x_1,x_2,x_3)=[Nx_1,Nx_2,Nx_3].
\end{eqnarray}

It follows from \eqref{eq:Nijenhuis1}, \eqref{eq:Nijenhuis2} and \eqref{eq:Nijenhuis3} that $N$ must satisfy the following condition
\begin{eqnarray}
\notag[N(x_1),N(x_2),N(x_3)]&=&N([Nx_1,Nx_2, x_3] + [Nx_1, x_2,Nx_3] + [x_1,Nx_2,Nx_3]) \\
\notag&&-N^2([Nx_1, x_2, x_3]+[x_1,Nx_2, x_3]+[x_1,x_2,Nx_3])\\
\label{eq:Nijenhuis4}&&+N^3([x_1,x_2,x_3]).
\end{eqnarray}

\begin{definition}(\cite{LSZB})
A linear operator $N:\frkg\to \frkg$ is called a \emph{Nijenhuis operator} if and only if \eqref{eq:Nijenhuis4} holds.
\end{definition}

We have seen that any second order trivial deformation produces a Nijenhuis operator.
Conversely, any Nijenhuis operator gives a second order trivial deformation as the following Theorem \ref{thm:Nijenhuis1} show.
The proof of the following Theorems \ref{thm:Nijenhuis1}, \ref{thm:Nijenhuis2} are by direct computations, see also \cite{LSZB,Zhang2}.

\begin{thm}\label{thm:Nijenhuis1}(\cite{LSZB})
Let $N$ be a Nijenhuis operator for $\frkg$. Then a  second order deformation of $\g$ can be obtained by putting
\begin{eqnarray*}
\omega_1(x_1,x_2,x_3)&=&[Nx_1,x_2,x_3]+[x_1,Nx_2,x_3]+[x_1,x_2,Nx_3]-N[x_1,x_2,x_3],\\
\omega_2(x_1,x_2,x_3)&=&[Nx_1,Nx_2,x_3]+[Nx_1,x_2,Nx_3]+[x_1,Nx_2,Nx_3]-N\omega_1(x_1,x_2,x_3).
\end{eqnarray*}
Furthermore, this deformation is a trivial one.
\end{thm}

\begin{thm}\label{thm:Nijenhuis2}(\cite{LSZB})
Let $N$ be a Nijenhuis operator. Then for any polynomial
$P(X)=\sum_{i=0}^n c_i X^i$, the operator $P(N)$ is also a Nijenhuis operator.
\end{thm}

\section{Extending structures for 3-Lie algebras}\label{sec:extending}
In this section, we study  extending structures for 3-Lie algebras.
Let $\g$ be a 3-Lie algebra and $E$ a vector space containing $\g$ as a subspace.
We are going to describe and classify all 3-Lie algebras structures on $E$ such that $\g$ is a subalgebra of
$E$. We show that associated to any extending structure of $\g$ by a complement space $V$, there is a unified product on the direct sum space  $E\cong\g\oplus V$.

\begin{definition}
 Let $(\frkg, [\cdot,\cdot,\cdot])$, $(E, [\cdot,\cdot,\cdot])$ be 3-Lie algebras and
$i:\frkg\to E$ be an inclusion map. The following short exact sequence (as vector spaces)
\begin{equation}\label{diagram:exact}
 \xymatrix{
   0  \ar[r]^{} & \frkg \ar[r]^{i} & E \ar[r]^{p} & V  \ar[r]^{} & 0 \\
  }
\end{equation}
is called  an  extending structures of $\frkg$ by $V$ if $E$  containing $\g$ as a subalgebra.
\end{definition}

\begin{definition}\label{lem:unifiedprod}
Let $\frkg$ be a 3-Lie algebra and $V$  a vector space.
An extending datum of $\g$ by $V$ is a system $\Omega(\g, V)$ consisting four maps
\begin{eqnarray*}
\triangleright:V\times V\times \g\to \g,\quad \triangleleft:V\times \g\times \g\to V, \\
\ppr:V\times \g\times \g\to \g,\quad \ppl:V\times V\times \g \to V,
\end{eqnarray*}
 and two totally skew-symmetric maps
 \begin{eqnarray*}
\omega:V\times V\times V\to \g, \quad \{\cdot,\cdot,\cdot\}:V\times V\times V\to V.
\end{eqnarray*}
Let $\Omega(\g,V)=(\trr, \trl,\ppr,\ppl,\omega,\{\cdot,\cdot,\cdot\})$ be an extending datum. Denote by $\g\natural V$ the direct sum vector space under the following bracket:
\begin{eqnarray}\label{eq:unifiedprod}
\notag&&[x_1 + u_1, x_2 + u_2, x_3 + u_3]\\
\notag&=&[x_1, x_2, x_3]+ (u_1, u_2)\trr x_3 +(u_2, u_3)\trr x_1+ (u_3, u_1)\trr x_2\\
\notag&&+u_3\ppr (x_1, x_2)+ u_1\ppr(x_2, x_3) +u_2\ppr (x_3, x_1)+\omega(u_1, u_2, u_3)\\
\notag&&+\{u_1, u_2, u_3\}+ u_1\trl(x_2, x_3)+ u_2\trl(x_3, x_1) + u_3 \trl(x_1, x_2) \\
&&+ (u_2, u_3)\ppl x_1+ (u_3, u_1)\ppl x_2+ (u_1, u_2)\ppl x_3,
\end{eqnarray}
where $x_i \in \frkg$ and $u_i \in V$.
Then $\frkg\natural V$ is called the unified product of $\frkg$ and $V$ if it is 3-Lie algebra with the above bracket.
\end{definition}

We will prove that there is a one-to-one correspondence between an  extending structures of $\frkg$ by $V$
and the unified product $\frkg\natural V$.

\begin{thm}\label{thm:unifiedprod}
Let $\frkg$ be a 3-Lie algebra, $V$  a vector space and $\Omega(\g,V)$ an extending datum.
Then $\g\natural\h$ is a unified product if and only if the following compatibility conditions hold:

\begin{eqnarray}\label{unifyI(1)-1}
&&v_3\trl[(x_1, x_2),(y_1, y_2)]=(v_3\trl(x_1, x_2))\trl(y_1, y_2)-(v_3\trl(y_1, y_2))\trl(x_1, x_2),
\end{eqnarray}
\begin{eqnarray}\label{unifyI(1)-2}
\notag && v_3\ppr[(x_1, x_2) ,(y_1, y_2)]\\
\notag&=&[x_1, x_2,v_3\ppr (y_1, y_2)]+(v_3\trl(y_1, y_2)) \ppr(x_1, x_2)\\
&&-[y_1, y_2, v_3\ppr (x_1, x_2)]-(v_3\trl(x_1, x_2)) \ppr(y_1, y_2),
\end{eqnarray}

\begin{eqnarray}\label{unifyI(2)-1}
\notag u_1\triangleleft (x_2, [y_1, y_2 ,y_3])
&=&(u_1\triangleleft (x_2, y_1))\trl (y_2, y_3)+(u_1\triangleleft (x_2, y_2))\trl (y_3, y_1)\\
&&+(u_1\triangleleft (x_2, y_3))\trl (y_1, y_2),
\end{eqnarray}
\begin{eqnarray}\label{unifyI(2)-2}
\notag &&u_1\ppr (x_2, [y_1, y_2 ,y_3])\\
\notag&=&[u_1\ppr (x_2, y_1),y_2, y_3] +(u_1\triangleleft (x_2, y_1))\ppr (y_2, y_3)\\
\notag&&+[u_1\ppr (x_2, y_2),y_3, y_1] +(u_1\triangleleft (x_2, y_2))\ppr (y_3, y_1)\\
&&+[u_1\ppr (x_2, y_3),y_1, y_2] +(u_1\triangleleft (x_2, y_3))\ppr (y_1, y_2),
\end{eqnarray}

\begin{eqnarray}\label{unifyII(1)-1}
\notag&&[x_1, x_2,(v_1,v_2)\triangleright y_3]+((v_1,v_2)\ppl y_3)\ppr (x_1, x_2)\\
\notag&=&v_2\ppr\big(v_1\ppr (x_1, x_2), y_3\big)+\big(v_1\trl (x_1, x_2), v_2\big)\trr y_3\\
\notag&&+v_1\ppr \big(v_2\ppr (x_1, x_2), y_3\big)+\big(v_1, v_2\trl(x_1, x_2)\big)\trr y_3\\
&&+(v_1, v_2)\trr[x_1, x_2, y_3],
\end{eqnarray}
\begin{eqnarray}\label{unifyII(1)-2}
\notag&&((v_1,v_2)\ppl y_3)\trl (x_1, x_2)\\
\notag&=&v_2\trl\big(v_1\ppr (x_1, x_2), y_3\big)+\big(v_1\trl (x_1, x_2), v_2\big)\ppl y_3\\
\notag&&+v_1\trl \big(v_2\ppr (x_1, x_2), y_3\big)+\big(v_1, v_2\trl(x_1, x_2)\big)\ppl y_3\\
&&+(v_1, v_2)\ppl[x_1, x_2, y_3],
\end{eqnarray}

\begin{eqnarray}
\notag&&u_1\ppr\big(x_2,v_1\ppr(y_2,y_3)\big)+\big(v_1,u_1\trl(x_2,y_3)\big)\trr y_2\\
\notag&&+[(u_1,v_1)\trr x_2,y_2,y_3]+((u_1,v_1)\ppl x_2)\ppr(y_2,y_3)\\
\notag&=&v_1\ppr\big(u_1\ppr(x_2,y_2),y_3\big)+\big(v_1, u_1\trl(x_2,y_2)\big)\trr y_3\\
&&+v_1\ppr\big(y_2,u_1\ppr(x_2,y_3)\big)+\big(u_1,v_1\trl(y_2,y_3)\big)\trr x_2,
\end{eqnarray}
\begin{eqnarray}
\notag&&u_1\trl\big(x_2,v_1\ppr(y_2,y_3)\big)+((u_1,v_1)\ppl x_2)\trl(y_2,y_3)\\
\notag&&+\big(v_1,u_1\trl(x_2,y_3)\big)\ppl y_2\\
\notag&=&\big(u_1,v_1\trl(y_2,y_3)\big)\ppl x_2+v_1\trl\big(u_1\ppr(x_2,y_2),y_3\big)\\
&&+\big(v_1, u_1\trl(x_2,y_2)\big)\ppl y_3+v_1\trl\big(y_2,u_1\ppr(x_2,y_3)\big),
\end{eqnarray}

\begin{eqnarray}
\notag&&(u_1, u_2)\triangleright[y_1, y_2 ,y_3]\\
\notag&=&[(u_1, u_2)\triangleright y_1, y_2 ,y_3]+((u_1, u_2)\ppl y_1)\ppr(y_2 ,y_3)\\
\notag&&+[y_1, (u_1, u_2)\triangleright y_2 ,y_3]+((u_1, u_2)\ppl y_2)\ppr(y_3 ,y_1)\\
&&+[y_1, y_2,(u_1, u_2)\triangleright y_3]+((u_1, u_2)\ppl y_3)\ppr(y_1 ,y_2),
\end{eqnarray}
\begin{eqnarray}
\notag&&(u_1, u_2)\ppl[y_1, y_2 ,y_3]\\
\notag&=&((u_1, u_2)\ppl y_1)\trl(y_2 ,y_3)+((u_1, u_2)\ppl y_2)\trl(y_3 ,y_1)\\
&&+((u_1, u_2)\ppl y_3)\trl(y_1 ,y_2),
\end{eqnarray}

\begin{eqnarray}
\notag&&\{v_1, v_2 ,v_3\}\ppr(x_1, x_2)+[x_1, x_2,\omega(v_1, v_2 ,v_3)]\\
\notag&=&(v_2, v_3)\trr (v_1\ppr(x_1, x_2))+\omega(v_1\trl(x_1, x_2), v_2, v_3)\\
\notag&&+(v_3, v_1)\trr (v_2\ppr(x_1, x_2))+\omega(v_2\trl(x_1, x_2), v_3, v_1)\\
&&+(v_1, v_2)\trr (v_3\ppr(x_1, x_2))+\omega(v_3\trl(x_1, x_2), v_1, v_2),
\end{eqnarray}
\begin{eqnarray}
\notag&&\{v_1, v_2 ,v_3\}\trl(x_1, x_2)\\
\notag&=&(v_2, v_3)\ppl (v_1\ppr(x_1, x_2))+\{v_1\trl(x_1, x_2), v_2, v_3\}\\
\notag&&+(v_3, v_1)\ppl (v_2\ppr(x_1, x_2))+\{v_2\trl(x_1, x_2), v_3, v_1\}\\
&&+(v_1, v_2)\ppl (v_3\ppr(x_1, x_2))+\{v_3\trl(x_1, x_2), v_1, v_2\},
\end{eqnarray}

\begin{eqnarray}
\notag&&u_2\ppr\big((v_1, v_2)\trr y_3, x_1\big)+\big(u_2,(v_1, v_2)\ppl y_3\big)\trr x_1\\
\notag&=&v_2\ppr\big(y_3,(u_2 ,v_1)\trr x_1\big)+\big((u_2 ,v_1)\ppl x_1, v_2\big)\trr y_3\\
\notag&&+v_1\ppr\big((u_2, v_2)\trr x_1, y_3\big)+\big(v_1, (u_2, v_2)\ppl x_1)\trr y_3\\
&&+(v_1, v_2)\trr\big(u_2\ppr(y_3,x_1)\big)+\omega(v_1, v_2,u_3\trl(y_3,x_1)),
\end{eqnarray}

\begin{eqnarray}
\notag&&+u_2\trl\big((v_1, v_2)\trr y_3, x_1\big)+\big(u_2,(v_1, v_2)\ppl y_3\big)\ppl x_1\\
\notag&=&v_2\trl\big(y_3,(u_2 ,v_1)\trr x_1\big)+\big((u_2 ,v_1)\ppl x_1, v_2\big)\ppl y_3\\
\notag&&+v_1\trl\big((u_2, v_2)\trr x_1, y_3\big)+\big(v_1, (u_2, v_2)\ppl x_1)\ppl y_3\\
&&+(v_1, v_2)\ppl \big(u_2\ppr(y_3,x_1)\big)+\{v_1, v_2,u_3\trl(y_3,x_1)\},
\end{eqnarray}

\begin{eqnarray}
\notag&&(u_1, u_2)\trr \big(v_1\ppr(y_2 ,y_3)\big)+\omega(u_1, u_2, v_1\trl(y_2 ,y_3))\\
\notag&=&[\omega(u_1, u_2 ,v_1), y_2, y_3]+\{u_1, u_2 ,v_1\}\ppr(y_2, y_3)\\
\notag&&+v_1\ppr\big((u_1, u_2)\trr y_2,  y_3\big)+\big(v_1, (u_1, u_2)\ppl y_2\big)\trr y_3\\
&&+v_1\ppr\big(y_2,(u_1, u_2)\trr y_3\big)+\big((u_1, u_2)\ppl y_3,v_1\big)\trr y_2,
\end{eqnarray}

\begin{eqnarray}
\notag&&(u_1, u_2)\ppl \big(v_1\ppr(y_2 ,y_3)\big)+\{u_1, u_2, v_1\trl(y_2 ,y_3)\}\\
\notag&=&\{u_1, u_2 ,v_1\}\trl(y_2, y_3)\\
\notag&&+v_1\trl\big((u_1, u_2)\trr y_2,  y_3\big)+\big(v_1, (u_1, u_2)\ppl y_2\big)\ppl y_3\\
&&+v_1\trl\big(y_2,(u_1, u_2)\trr y_3\big)+\big((u_1, u_2)\ppl y_3,v_1\big)\ppl y_2,
\end{eqnarray}

\begin{eqnarray}
\notag&&\{(u_1, u_2),(v_1, v_2)\}\trr y_3\\
\notag&=&(u_1, u_2)\trr((v_1, v_2)\trr y_3)-(v_1, v_2)\trr((u_1, u_2)\trr y_3)\\
\notag&&+\omega(u_1, u_2,(v_1, v_2)\ppl y_3)-\omega(v_1, v_2,(u_1, u_2)\ppl y_3)\\
&&+v_2\ppr(\omega(u_1, u_2 ,v_1), y_3)-v_1\ppr(\omega(u_1, u_2 ,v_2), y_3),
\end{eqnarray}
\begin{eqnarray}
\notag&&\{(u_1, u_2), (v_1, v_2)\}\ppl y_3\\
\notag&=&(u_1, u_2)\ppl((v_1, v_2)\trr y_3)+[u_1, u_2,(v_1, v_2)\ppl y_3]\\
\notag&&-(v_1, v_2)\ppl((u_1, u_2)\trr y_3)-[v_1, v_2,(u_1, u_2)\ppl y_3]\\
&&+v_2\trl(\omega(u_1, u_2 ,v_1), y_3)-v_1\trl(\omega(u_1, u_2 ,v_2), y_3),
\end{eqnarray}

\begin{eqnarray}
\notag&&(u_2, [v_1, v_2 ,v_3])\trr x_1+u_2\ppr(\omega(v_1, v_2 ,v_3),x_1)\\
\notag&=&(v_2, v_3)\trr((u_2, v_1)\trr x_1)+\omega((u_2, v_1)\ppl x_1,v_2, v_3)\\
\notag&&+ (v_3, v_1)\trr((u_2, v_2)\trr x_1)+\omega((u_2, v_2)\ppl x_1,v_3, v_1)\\
&&+(v_1, v_2)\trr((u_2, v_3)\trr x_1)+\omega((u_2, v_3)\ppl x_1,v_1, v_2),
\end{eqnarray}
\begin{eqnarray}
\notag&& (u_2, [v_1, v_2 ,v_3])\ppl x_1+u_2\trl(\omega(v_1, v_2 ,v_3),x_1)\\
\notag&=&(v_2, v_3)\trr((u_2, v_1)\ppl x_1) +\{(u_2, v_1)\ppl x_1,v_2, v_3\}\\
\notag&&+ (v_3, v_1)\trr((u_2, v_2)\ppl x_1) +\{(u_2, v_2)\ppl x_1,v_3, v_1\}\\
&&+ (v_1, v_2)\trr((u_2, v_3)\ppl x_1) +\{(u_2, v_3)\ppl x_1,v_1, v_2\},
\end{eqnarray}

\begin{eqnarray}
\notag&&\{u_1, u_2,\{v_1, v_2 ,v_3\}\}-\{\{u_1, u_2,v_1\}, v_2 ,v_3\}\}\\
\notag&&-\{v_1,\{u_1, u_2, v_2\} ,v_3\}\}-\{v_1, v_2,\{u_1, u_2 ,v_3\}\}\\
\notag&&+(u_1, u_2)\ppl\omega(v_1, v_2 ,v_3)-(v_2, v_3)\ppl\omega(u_1, u_2,v_1)\\
&&-(v_3, v_1)\ppl\omega(u_1, u_2,v_2)-(v_1, v_2)\ppl\omega(u_1, u_2 ,v_3)=0,
\end{eqnarray}
\begin{eqnarray}
\notag&&\omega(u_1, u_2,\{v_1, v_2 ,v_3\})-\omega(\{u_1, u_2,v_1\} v_2 ,v_3)\\
\notag&&-\omega(v_2, \{u_1, u_2,v_2\} ,v_3)-\omega(v_1, v_2,\{u_1, u_2 ,v_3\})\\
\notag&&+(u_1, u_2)\trr\omega(v_1, v_2 ,v_3)-(v_2, v_3)\trr\omega(u_1, u_2,v_1)\\
&&-(v_3, v_1)\trr\omega(u_1, u_2,v_2)-(v_1, v_2)\trr\omega(u_1, u_2 ,v_3)=0.
\end{eqnarray}
In the above conditions, \eqref{unifyI(1)-1} and \eqref{unifyI(2)-1} means that $(V,\trl)$ is a representation of $\g$.
\end{thm}

\pf
Assume $\g\natural\h$ is a unified product, then we have
\begin{eqnarray*}
&&[x_1 + u_1, x_2 + u_2,[y_1 + v_1, y_2 + v_2, y_3 + v_3]]\\
&=&[[x_1 + u_1, x_2 + u_2,y_1 + v_1], y_2 + v_2, y_3 + v_3]\\
&&+[y_1 + v_1, [x_1 + u_1, x_2 + u_2, y_2 + v_2], y_3 + v_3]\\
&&+[y_1 + v_1, y_2 + v_2,[x_1 + u_1, x_2 + u_2, y_3 + v_3]].
\end{eqnarray*}
By direct computations of the fundamental identity for different combinations
of elements, assuming some of $x_i$ and $u_i$ are equal to zero,  we can check the above conditions must hold in the above five cases.
For example,  by letting $u_1,u_2, y_1,, y_2, v_3$ to be zero, we have the equality
\begin{eqnarray*}
&&[x_1, x_2,[v_1, v_2 ,y_3]]=[[x_1, x_2 ,v_1], v_2, y_3]+[v_1, [x_1, x_2, v_2],  y_3]+[v_1, v_2,[x_1, x_2, y_3]].
\end{eqnarray*}
By definition of unified product \eqref{eq:unifiedprod}, we get
\begin{eqnarray*}
[x_1, x_2,[v_1, v_2 ,y_3]]&=&[x_1, x_2,(v_1,v_2)\triangleright y_3]+((v_1,v_2)\ppl y_3)\ppr (x_1, x_2)\\
&&+((v_1,v_2)\ppl y_3)\trl (x_1, x_2),\\
{[[x_1, x_2 ,v_1], v_2, y_3]}&=&[v_1\ppr (x_1, x_2), v_2, y_3]+[v_1\trl (x_1, x_2), v_2, y_3]\\
&=&v_2\ppr(v_1\ppr (x_1, x_2), y_3)+v_2\trl(v_1\ppr (x_1, x_2), y_3)\\
&&+(v_1\trl (x_1, x_2), v_2)\trr y_3+(v_1\trl (x_1, x_2), v_2)\ppl y_3,\\
{[v_1, [x_1, x_2, v_2], y_3]}&=&[v_1, v_2\ppr (x_1, x_2), y_3]+[v_1, v_2\trl(x_1, x_2), y_3]\\
&=&v_1\ppr (v_2\ppr (x_1, x_2), y_3)+v_1\trl (v_2\ppr (x_1, x_2), y_3)\\
&&+(v_1, v_2\trl(x_1, x_2))\trr y_3+(v_1, v_2\trl(x_1, x_2))\ppl y_3,\\
{[v_1, v_2,[x_1, x_2, y_3]]}&=&(v_1, v_2)\trr[x_1, x_2, y_3]+(v_1, v_2)\ppl([x_1, x_2, y_3]).
\end{eqnarray*}
Thus we obtain
\begin{eqnarray}
\notag&&[x_1, x_2,(v_1,v_2)\triangleright y_3]+((v_1,v_2)\ppl y_3)\ppr (x_1, x_2)\\
\notag&=&v_2\ppr(v_1\ppr (x_1, x_2), y_3)+(v_1\trl (x_1, x_2), v_2)\trr y_3\\
\notag&&+v_1\ppr (v_2\ppr (x_1, x_2), y_3)+(v_1, v_2\trl(x_1, x_2))\trr y_3\\
&&+(v_1, v_2)\trr([x_1, x_2, y_3]),
\end{eqnarray}
and
\begin{eqnarray}
\notag&&((v_1,v_2)\ppl y_3)\trl (x_1, x_2)\\
\notag&=&v_2\trl(v_1\ppr (x_1, x_2), y_3)+(v_1\trl (x_1, x_2), v_2)\ppl y_3\\
\notag&&+v_1\trl (v_2\ppr (x_1, x_2), y_3)+(v_1, v_2\trl(x_1, x_2))\ppl y_3\\
&&+(v_1, v_2)\ppl([x_1, x_2, y_3]),
\end{eqnarray}
which are conditions \eqref{unifyII(1)-1} and \eqref{unifyII(1)-2} respectively.
The other conditions can obtained similarly.

Conversely, if the above conditions hold, it is straightforward to see that $\g\natural\h$  is a 3-Lie
algebra.
\qed\\

\begin{ex}
Let $\mathfrak{g}$ be the  three dimensional 3-Lie algebra  with zero bracket under the basis $\{x_1, x_2, x_3\}$,
 and  $V$ be a  five dimension vector space with basis $\{u_1, u_2, u_3, u_4, u_5\}$. We define the following extending datum:
$$
\begin{gathered}
\left[u_{2}, u_{3}, u_{4}\right]=u_{1},\quad \left[u_{1}, u_{3}, u_{4}\right]=u_{2}, \\
u_{5}\ppr (x_{1}, x_{2})=  x_{2}+c_{2} x_{3},\quad u_{5}\ppr (x_{1}, x_{3})=c_1 x_{2}-  x_{3},\\
\omega(u_{2}, u_{3}, u_{5})=-x_{3},\quad \omega(u_{3}, u_{4}, u_{5})=-x_{1},\quad \omega(u_{2}, u_{4}, u_{5})=x_{2}, \\
(u_{1}, u_{4})\trr x_{1}=x_{2},\quad (u_{1}, u_{3})\trr x_{1}=-x_{3},\quad (u_{3}, u_{4})\trr x_{1}=-x_{1}, \\
(u_{3}, u_{5})\ppl x_{1}=-  u_{3}-c_1 u_{4},\quad (u_{4}, u_{5})\ppl x_{1}= -c_{2} u_{3}+ u_{4}, \quad (u_{3}, u_{5})\ppl x_{2}= u_{2},\\
 (u_{4}, u_{5})\ppl x_{2}= c_2u_{2},\quad (u_{3}, u_{5})\ppl  x_{3}=c_1 u_{2},\quad (u_{4}, u_{5})\ppl  x_{3}=-  u_{2},\\
u_{3}\trl ( x_{1}, x_{2})=-  u_{1},\quad u_{4}\trl ( x_{1}, x_{2})=-  c_2 u_{1},\\
u_{3}\trl (x_{1}, x_{3})=-c_1   u_{1},\quad u_{4}\trl (x_{1}, x_{3})=  u_{1}.
\end{gathered}
$$
where $c_1, c_2$ are arbitrary parameter. Then one checks that this extending datum satisfying the conditions in the above  Theorem \ref{thm:unifiedprod}.

Thus we obtain an eight dimensional 3-Lie algebra $E=\mathfrak{g}\natural V$ defined with respect to a basis
$\{x_1, x_2, x_3,u_1, u_2, u_3,  u_4, u_5\}$ by the skew-symmetric bracket
$$
\begin{gathered}
{\left[u_{2}, u_{3}, u_{4}\right]=u_{1},\quad \left[u_{1}, u_{3}, u_{4}\right]=u_{2},} \\
{\left[u_{5}, x_{1}, x_{2}\right]=  x_{2}+c_{2} x_{3},\quad \left[u_{5}, x_{1}, x_{3}\right]=c_{1} x_{2}-  x_{3},} \\
{\left[u_{2}, u_{3}, u_{5}\right]=-x_{3},\quad \left[u_{3}, u_{4}, u_{5}\right]=-x_{1},\quad \left[u_{2}, u_{4}, u_{5}\right]=x_{2},} \\
{\left[u_{1}, u_{4}, x_{1}\right]=x_{2},\quad \left[u_{1}, u_{3}, x_{1}\right]=-x_{3},\quad \left[u_{3}, u_{4}, x_{1}\right]=-x_{1},} \\
{\left[u_{3}, u_{5}, x_{1}\right]=-  u_{3}-c_{1} u_{4},\quad \left[u_{4}, u_{5}, x_{1}\right]=-c_{2} u_{3}+  u_{4},} \\
{\left[u_{3}, u_{5}, x_{2}\right]=  u_{2},\quad\left[u_{4}, u_{5}, x_{2}\right]=c_{2} u_{2},} \\
{\left[u_{3}, u_{5}, x_{3}\right]=c_{1} u_{2},\quad\left[u_{4}, u_{5}, x_{3}\right]=-  u_{2},} \\
{\left[u_{3}, x_{1}, x_{2}\right]=-  u_{1},\quad \left[u_{4}, x_{1}, x_{2}\right]=-c_{2} u_{1},} \\
{\left[u_{3}, x_{1}, x_{3}\right]=-c_{1}   u_{1},\quad \left[u_{4}, x_{1}, x_{3}\right]=  u_{1}}.
\end{gathered}
$$
\end{ex}

By a suitable change of basis, this eight dimensional 3-Lie algebra  is isomorphic to the 3-Lie algebra given in \cite[Theorem 4.10(3)]{DBL} by using the $r$-matrix method, but here we have reconstructed it  by using a more direct approach. More examples of unified product for 3-Lie algebras with both non-zero brackets  on $\frkg$ and $V$ will be given in the next section.

Given an  extending structure $\Omega(\frkg,V)$, it is obvious that $\frkg$ can be seen as a 3-Lie subalgebra of $\frkg\natural V$. Conversely, we will prove that any 3-Lie algebra structure on a vector space $E$ containing $\frkg$ as a subalgebra is isomorphic to a unified product.

\begin{thm}\label{pt1}
Let $(\frkg,[\cdot,\cdot,\cdot]$ and $(E,[\cdot,\cdot,\cdot])$ be two 3-Lie algebras  such that $E$ containing $\frkg$ as a subalgebra of $E$. Then, there exists an extending structure $\Omega(\frkg,V)$ of $\frkg$ by a subspace $V$ of $E$ and an isomorphism of 3-Lie algebras $E\cong \frkg\natural V$ which stabilizes $\frkg$ and co-stabilizes $V$.
\end{thm}

\begin{proof}
Note that there is a natural linear map $p: E\rightarrow A$ such that $p(x)=x$ for all $x\in \g$.
Set $V=\text{Ker}(p)$ which is a complement of $\frkg$ in $E$. Then, we define the extending datum $\Omega(\frkg,V)$ of $\frkg$ by a subspace $V$ of $E$ as follows:
\begin{eqnarray*}
\triangleright:V\times V\times \g\to \g,&&\quad (u_1, u_2)\triangleright x_3=p([u_1, u_2, x_3]), \\
\ppr:V\times \g\times \g\to \g,&&\quad u_1\ppr (x_2, x_3)=p([u_1, x_2, x_3]),\\
\omega:V\times V\times V\to \g, &&\quad \omega(u_1, u_2, u_3)=p([u_1, u_2, u_3]),\\
 \triangleleft:V\times \g\times \g\to V, &&\quad u_1\triangleleft (x_2, x_3)=[u_1, x_2, x_3]-p([u_1, x_2, x_3]),\\
 \ppl:V\times V\times \g \to V,&&\quad (u_1, u_2)\ppl x_3=[u_1, u_2, x_3]-p([u_1, u_2, x_3]),\\
 \{\cdot,\cdot,\cdot\}:V\times V\times V\to V,&&\quad \{u_1, u_2, u_3\}=[u_1, u_2, u_3]-p([u_1, u_2, u_3]),
\end{eqnarray*}
for any $x_i \in \frkg, u_i \in V$. It is easy to see that $\varphi: \frkg\times V\rightarrow E$ defined as
$\varphi(x,u)=x+u$ is a linear isomorphism, whose inverse is as follows: $\varphi^{-1}(e):=(p(e),e-p(e))$ for all $e\in E$.
Next, we should prove that $\Omega(\frkg,V)$ is an  extending structure of $\frkg$ by $V$ and $\varphi: A \natural V \rightarrow E$ is an isomorphism of 3-Lie algebras that stabilizes $\frkg$ and co-stabilizes $V$. In fact, if $\varphi: \frkg\times V\rightarrow E$ is an isomorphism of 3-Lie algebras, there exists
a unique 3-bracket  given by
\begin{eqnarray}\label{ee1}
[(x_1, u_1), (x_2, u_2), (x_3, u_3)]=\varphi^{-1}[\varphi(x_1, u_1), \varphi(x_2, u_2), \varphi(x_3, u_3)]
\end{eqnarray}
We are going to proof that the 3-bracket defined by (\ref{ee1}) is just the one given in  the above extending system $\Omega(\frkg,V)$.
Indeed, for any $x_i \in \frkg, u_i \in V$,  we have
\begin{eqnarray*}
&&[(x_1, u_1), (x_2, u_2), (x_3, u_3)]\\
&=&\varphi^{-1}[\varphi(x_1, u_1), \varphi(x_2, u_2), \varphi(x_3, u_3)]\\
&=&\varphi^{-1}[x_1 + u_1, x_2 + u_2, x_3 + u_3]\\
&=&\varphi^{-1}\big([x_1, x_2, x_3]+{[x_1, x_2, u_3]}+{[x_1, u_2, u_3]}\\
&&+{[u_1, x_2, x_3]}+{[u_1, u_2, x_3]}+{[u_1, u_2, u_3]}\big)\\
&=&\Big([x_1, x_2, x_3]+p([u_1, u_2, x_3])\\
&&+p([u_1, x_2, x_3])+c.p.+p([u_1, u_2, u_3])+c.p.,\\
&&\quad[u_1, u_2, u_3]+[u_1, u_2, x_3]+c.p.+[u_1, x_2, x_3]+c.p.\\
&&\quad-p([u_1, u_2, u_3])-p([u_1, x_2, x_3])-c.p.-p([u_1, u_2, x_3]-c.p.\Big)\\
&=&\Big([x_1, x_2, x_3]+ (u_1, u_2)\trr x_3 +(u_2, u_3)\trr x_1+ (u_3, u_1)\trr x_2\\
&&+u_3\ppr (x_1, x_2)+ u_1\ppr(x_2, x_3) +u_2\ppr (x_3, x_1)+\omega(u_1, u_2, u_3),\\
&&\quad\{u_1, u_2, u_3\}+ u_1\trl(x_2, x_3)+ u_2\trl(x_3, x_1) + u_3 \trl(x_1, x_2) \\
&&\quad+ (u_2, u_3)\ppl x_1+ (u_3, u_1)\ppl x_2+ (u_1, u_2)\ppl x_3\Big).
\end{eqnarray*}
This bracket coincides with the bracket in (\ref{eq:unifiedprod})   and thus $\varphi: \frkg\natural V\rightarrow E$ is an isomorphism of 3-Lie algebras.
Moreover, the following diagram is commutative
$$\xymatrix{ {\frkg}\ar[d]^{id}\ar[r]^{i}& {\frkg\natural V}\ar[d]^{\varphi}\ar[r]^{q} & {V}\ar[d]^{id} \\
{\frkg}\ar[r]^{i}&{E}\ar[r]^{\pi} &{V} }$$
where $q:\frkg\natural V\rightarrow V$ and $\pi: E\rightarrow V$ are the natural projections.
The proof is finished.
\end{proof}

\begin{definition}
Let $\frkg$ be a 3-Lie algebra, $E$ a vector space such that $\frkg$ is a subspace
of $E$ and $V$ a complement of $\frkg$ in $E$. For a linear map $\varphi: E\rightarrow E$, the following
diagram is considered:
$$\xymatrix{ {\frkg}\ar[d]^{id}\ar[r]^{i}& {E}\ar[d]^{\varphi}\ar[r]^{\pi} & {V}\ar[d]^{id} \\
{\frkg}\ar[r]^{i}&{E}\ar[r]^{\pi} &{V} }$$
where $\pi: E\rightarrow V$ is the natural projection of $E=\frkg\oplus V$ onto $V$ and
$i: \frkg\rightarrow E$ is the inclusion map. We say that
$\varphi: E\rightarrow E$ \emph{stabilizes} $\frkg$ (resp. \emph{co-stabilizes} $V$) if the left square
(resp. the right square) of the  above diagram is commutative.

Let $[\cdot,\cdot,\cdot]$ and $[\cdot,\cdot,\cdot]'$ be two 3-Lie  structures on $E$ both containing $\frkg$ as a 3-Lie   subalgebra.
If there exists a 3-Lie algebra isomorphism $\varphi: (E, [\cdot,\cdot,\cdot])\rightarrow (E,[\cdot,\cdot,\cdot]')$ which stabilizes
$\frkg$, then $(E, [\cdot,\cdot,\cdot])$ and $(E,[\cdot,\cdot,\cdot]')$ are called \emph{equivalent}, which is denoted by
$(E, [\cdot,\cdot,\cdot])\equiv (E,[\cdot,\cdot,\cdot]')$.
If there exists a 3-Lie algebra isomorphism $\varphi: (E, [\cdot,\cdot,\cdot])\rightarrow (E,[\cdot,\cdot,\cdot]')$ which stabilizes
$\frkg$ and co-stabilizes $V$,   then $(E, [\cdot,\cdot,\cdot])$ and $(E,[\cdot,\cdot,\cdot]')$  are called \emph{cohomologous}, which is denoted by
$(E, [\cdot,\cdot,\cdot])\approx (E,[\cdot,\cdot,\cdot]')$.
\end{definition}

Obviously, $\equiv $ and $\approx $ are equivalence relations on the set of all 3-Lie algebra structures on
$E$ containing $\frkg$ as a 3-Lie subalgebra. Denote by $\text{Extd}(E,\frkg)$ (resp. $\text{Extd}'(E,\frkg)$) the set
of all equivalence classes via $\equiv $ (resp. $\approx $). Thus, $\text{Extd}(E,\frkg)$ is the classifying object of the extending structures problem and $\text{Extd}'(E,\frkg)$ provides a classification of the extending structures problem from the point of the view of the extension problem. In addition, it is easy to see that there exists a canonical projection $\text{Extd}(E,\frkg)\twoheadrightarrow \text{Extd}'(E,\frkg)$.

Next, by Theorem \ref{pt1}, for classifying all 3-Lie algebra structures on $E$ containing $\frkg$ as a subalgebra,
we only need to classify all unified products $\frkg\natural V$ associated to all 3-Lie algebra structures
$\Omega(\frkg,V) $ for a given complement $V$ of $\frkg$ in $E$.

\begin{lem}\label{mainlemma}
Let $\Omega(\g, V) = (\trr, \trl,\ppr, \ppl, \omega, \{\cdot,\cdot,\cdot\})$ and
$\Omega'(\g, V) = (\trr', \trl',\ppr', \ppl', \omega', \{\cdot,\cdot,\cdot\}')$  be two  extending structures of $\g$ through $V$
and $\g\natural V , \g\natural' V$ the associated unified products. Then there exists a bijection between the set of all homomorphisms
of 3-Lie algebra  $\psi : \g\natural V \to \g\natural'V$ which stabilizes $\g$ and the set of pairs $(r, \nu)$,
where $r : V \to \g, \nu: V\to V$ are two linear maps satisfying the following compatibility
conditions for all $x_i \in \g, u_i\in V$:

(M1) $$\nu(u_1\trl(x_2, x_3))=\nu(u_1)\trl' (x_2, x_3),$$

(M2) $$r(u_1\trl(x_2, x_3))+u_1\ppr (x_2, x_3)=[x_2, x_3, r(u_1)]+\nu(u_1)\ppr' (x_2, x_3),$$

(M3)
\begin{eqnarray*}
&&\nu((u_1, u_2)\ppl x_3)\\
&=&\big(\nu(u_1), \nu(u_2)\big)\ppl' x_3-\nu(u_2)\trl' (r(u_1), x_3)+\nu(u_1)\trl' (r(u_2), x_3),
\end{eqnarray*}

(M4)
\begin{eqnarray*}
&&r((u_1, u_2)\ppl x_3)+(u_1, u_2)\triangleright x_3\\
&=&[r(u_1), r(u_2), x_3]+(\nu(u_1), \nu(u_2))\trr' x_3\\
&&-\nu(u_2)\ppr' (r(u_1), x_3)+\nu(u_1)\ppr' (r(u_2), x_3),
\end{eqnarray*}

(M5)
\begin{eqnarray*}
&&\nu(\{u_1, u_2, u_3\})\\
&=&\{\nu(u_1), \nu(u_2), \nu(u_3)\}'+(\nu(u_1), \nu(u_2))\ppl' r(u_3)+c.p.\\
&&+\nu(u_1)\trl'(r(u_2), r(u_3))+c.p.,
\end{eqnarray*}

(M6)
\begin{eqnarray*}
&&r(\{u_1, u_2, u_3\})+\omega(u_1, u_2, u_3)\\
&=&[r(u_1), r(u_2), r(u_3)]+\omega'(\nu(u_1), \nu(u_2), \nu(u_3))\\
&&+(\nu(u_1), \nu(u_2))\trr' r(u_3)+c.p.+\nu(u_1)\ppr'(r(u_2), r(u_3))+c.p.
\end{eqnarray*}
where c.p. means cyclic permutations with respect to elements $u_1, u_2$ and $u_3$.

Under the above bijection the homomorphism of 3-Lie algebras $\psi: \g\natural V\to \g\natural' V$
corresponding to $(r, \nu)$ is given by:
$$\psi(x, u) = (x + r(u), \nu(u)).$$
Moreover, $\psi$ is an isomorphism if and only if $\nu : V \to V$ is an isomorphism and
$\psi$ co-stabilizes V if and only if $\nu = id_V$.
\end{lem}

\pf
A linear map $\psi : \g\natural V \to \g \natural'V$  which stabilizes $\g$ is uniquely determined by
two linear maps $r : V \to \g,\ \nu: V\to V$ such that $\psi(x, u) = (x + r(u), \nu(u))$, for all $x \in \g$ and $u \in V$.
We rewrite the map $\psi$ as an equivalent forms as
\begin{eqnarray*}
\psi(x+u)&=&\psi(x)+\psi(u)\\
&=&x+r(u)+\nu(u).
\end{eqnarray*}
Now we prove that $\psi$ is a homomorphism if and only if (M1)--(M6) hold. By the equality
$$\psi([x_1 + u_1, x_2 + u_2, x_3 + u_3])=[\psi(x_1 + u_1), \psi(x_2 + u_2), \psi(x_3 + u_3)]', \qquad (*)$$
we have the following three cases.

Case I:
\begin{eqnarray*}
\psi([u_1,x_2, x_3])&=&\psi(u_1\ppr (x_2, x_3))+\psi(u_1\trl(x_2, x_3)),\\
&=&u_1\ppr (x_2, x_3)+r(u_1\trl(x_2, x_3))+\nu(u_1\trl(x_2, x_3)),
\end{eqnarray*}
\begin{eqnarray*}
{[\psi(u_1),\psi(x_2), \psi(x_3)]'}&=&[r(u_1),x_2, x_3]+[\nu(u_1),x_2, x_3]'\\
&=&[r(u_1),x_2, x_3]+\nu(u_1)\ppr' (x_2, x_3)+(\nu(u_1))\trl' (x_2, x_3).
\end{eqnarray*}
Thus we obtain (M1) and (M2).

Case II:
\begin{eqnarray*}
\psi([u_1, u_2, x_3])&=&\psi((u_1, u_2)\triangleright x_3)+\psi((u_1, u_2)\ppl x_3),\\
&=&(u_1, u_2)\triangleright x_3+r((u_1, u_2)\ppl x_3)+\nu((u_1, u_2)\ppl x_3),
\end{eqnarray*}
\begin{eqnarray*}
{[\psi(u_1), \psi(u_2), \psi(x_3)]'}&=&[r(u_1)+\nu(u_1), r(u_2)+\nu(u_2), x_3]'\\
&=&[r(u_1), r(u_2), x_3]+[\nu(u_1), \nu(u_2), x_3]'\\
&&+[r(u_1), \nu(u_2), x_3]'+[\nu(u_1), r(u_2), x_3]'\\
&=&[r(u_1), r(u_2), x_3]\\
&&+(\nu(u_1), \nu(u_2))\trr' x_3+(\nu(u_1), \nu(u_2))\ppl' x_3\\
&&-\nu(u_2)\ppr' (r(u_1), x_3)-\nu(u_2)\trl' (r(u_1), x_3)\\
&&+\nu(u_1)\ppr' (r(u_2), x_3)+\nu(u_1)\trl' (r(u_2), x_3).
\end{eqnarray*}
Thus we obtain (M3) and (M4).

Case III:
\begin{eqnarray*}
\psi([u_1, u_2, u_3])&=&\psi(\{u_1, u_2, u_3\})+\psi(\omega(u_1, u_2, u_3)),\\
&=&r(\{u_1, u_2, u_3\})+\nu(\{u_1, u_2, u_3\})+\omega(u_1, u_2, u_3),
\end{eqnarray*}
\begin{eqnarray*}
{[\psi(u_1), \psi(u_2), \psi(u_3)]'}&=&[r(u_1)+\nu(u_1), r(u_2)+\nu(u_2), r(u_3)+\nu(u_3)]'\\
&=&[r(u_1), r(u_2), r(u_3)]\\
&&+\nu(u_2)\ppr'(r(u_3), r(u_1))+\nu(u_2)\trl'(r(u_3), r(u_2)) \\
&&+\nu(u_1)\ppr'(r(u_2), r(u_3))+\nu(u_1)\trl'(r(u_2), r(u_3)) \\
&&+\nu(u_3)\ppr'(r(u_1), r(u_2))+\nu(u_3)\trl'(r(u_1), r(u_2)) \\
&&+(\nu(u_3), \nu(u_1))\trr' r(u_2)+(\nu(u_3), \nu(u_2))\ppl' r(u_2)\\
&&+(\nu(u_2), \nu(u_3))\trr' r(u_1)+(\nu(u_2), \nu(u_3))\ppl' r(u_1)\\
&&+(\nu(u_1), \nu(u_2))\trr' r(u_3)+(\nu(u_1), \nu(u_2))\ppl' r(u_3)\\
&&+\{\nu(u_1), \nu(u_2), \nu(u_3)\}'+\omega'(\nu(u_1), \nu(u_2), \nu(u_3)).
\end{eqnarray*}
Thus we obtain (M5) and (M6).
Conversely, if conditions (M1)--(M6) hold,  then by direct computations of both sides of $(*)$,
one check that the left hand side of $(*)$ is equal to the right hand side.

Assume that $\nu: V \to V$ is bijective. Then $\psi$ is
an isomorphism of 3-Lie algebras with the inverse given by $
\psi^{-1}(x, u) = \bigl(x - r(\nu^{-1}(u)),
\nu^{-1}(y)\bigl)$, for all $x \in \mathfrak{g}$ and $u \in V$.
Conversely, assume that $\psi$ is bijective. It follows
easily that $\nu$ is surjective. Now we prove that $v$
is injective. Let $u \in V$ such that $\nu(u) = 0$. We have
 $\psi(0, 0) = (0, 0) = (0, \nu(u)) = \psi(-
r(u), u)$, and hence we obtain $u = 0$, i.e. $\nu$ is a bijection.
Finally, it is straightforward to see that $\psi$ co-stabilizes $V$ if and only
if $\nu= id_V$ and the proof is now finished.
\qed

\begin{definition}
Two extending structures  $\Omega(\g, V) = (\trr, \trl,\ppr, \ppl, \omega, \{\cdot,\cdot,\cdot\})$ and
$\Omega'(\g, V) = (\trr', \trl',\ppr', \ppl', \omega', \{\cdot,\cdot,\cdot\}')$  are called  equivalent  if there exists  a pair $(r,\nu)$  of linear maps, where $r : V \to \g$ and  $\nu: V\to V$ is an isomorphism satisfying the following  conditions for all $x_i \in \g, u_i\in V$:
\begin{eqnarray}
u_1\trl(x_2, x_3)=\nu^{-1}\big(\nu(u_1)\trl' (x_2, x_3)\big),
\end{eqnarray}
\begin{eqnarray}
u_1\ppr (x_2, x_3)=[x_2, x_3, r(u_1)]+\nu(u_1)\ppr' (x_2, x_3)-r\circ\nu^{-1}\Big(\nu(u_1)\trl' (x_2, x_3)\Big),
\end{eqnarray}
\begin{eqnarray}
(u_1, u_2)\ppl x_3&=&\nu^{-1}\big((\nu(u_1), \nu(u_2))\ppl' x_3-\nu(u_2)\trl' (r(u_1), x_3)\\
&&+\nu(u_1)\trl' (r(u_2), x_3)\big),
\end{eqnarray}
\begin{eqnarray}
\notag(u_1, u_2)\triangleright x_3&=&[r(u_1), r(u_2), x_3]+(\nu(u_1), \nu(u_2))\trr' x_3-\nu(u_2)\ppr' (r(u_1), x_3)\\
\notag&&+\nu(u_1)\ppr' (r(u_2), x_3)-r\circ\nu^{-1}\Big((\nu(u_1), \nu(u_2))\ppl' x_3\\
&&-\nu(u_2)\trl' (r(u_1), x_3)+\nu(u_1)\trl' (r(u_2), x_3)\Big),
\end{eqnarray}
\begin{eqnarray}
\notag\{u_1, u_2, u_3\}&=&\nu^{-1}\Big(\{\nu(u_1), \nu(u_2), \nu(u_3)\}'+(\nu(u_1), \nu(u_2))\ppl' r(u_3)+c.p.\\
&&+\nu(u_1)\trl'(r(u_2), r(u_3))+c.p.\Big),
\end{eqnarray}
\begin{eqnarray}
\notag\omega(u_1, u_2, u_3)&=&[r(u_1), r(u_2), r(u_3)]+\omega'(\nu(u_1), \nu(u_2), \nu(u_3))\\
\notag&&+(\nu(u_1), \nu(u_2))\trr' r(u_3)+c.p.+\nu(u_1)\ppr'(r(u_2), r(u_3))+c.p.\\
\notag&&-r\circ\nu^{-1}\Big(\{\nu(u_1), \nu(u_2), \nu(u_3)\}'+(\nu(u_1), \nu(u_2))\ppl' r(u_3)+c.p.\\
&&+\nu(u_1)\trl'(r(u_2), r(u_3))+c.p.\Big),
\end{eqnarray}
where c.p. means cyclic permutations with respect to elements $u_1, u_2$ and $u_3$.
\end{definition}

\begin{definition}
Two extending structures  $\Omega(\g, V) = (\trr, \trl,\ppr, \ppl, \omega, \{\cdot,\cdot,\cdot\})$ and
$\Omega'(\g, V) = (\trr', \trl',\ppr', \ppl', \omega', \{\cdot,\cdot,\cdot\}')$  are called  cohomologous if $\trl=\trl', \ppl=\ppl'$ and there exists  a  linear maps $r : V \to \g$  satisfying the following  conditions for all $x_i \in \g, u_i\in V$:
\begin{eqnarray}
u_1\ppr (x_2, x_3)=u_1\ppr' (x_2, x_3)+[x_2, x_3, r(u_1)] -r(u_1\trl(x_2, x_3)),
\end{eqnarray}
\begin{eqnarray}
u_2\trl' (r(u_1), x_3)=u_1\trl' (r(u_2), x_3),
\end{eqnarray}
\begin{eqnarray}
\notag(u_1, u_2)\triangleright x_3&=&(u_1, u_2)\trr' x_3+[r(u_1), r(u_2), x_3]-r((u_1, u_2)\ppl' x_3)\\
&&- u_2\ppr' (r(u_1), x_3)+ u_1\ppr' (r(u_2), x_3),
\end{eqnarray}
\begin{eqnarray}
\notag\{u_1, u_2, u_3\}&=& \{u_1, u_2, u_3\}'+(u_1, u_2)\ppl' r(u_3)+c.p.\\
&&+ u_1\trl'(r(u_2), r(u_3))+c.p.,
\end{eqnarray}
\begin{eqnarray}
\notag\omega(u_1, u_2, u_3)&=&\omega'(u_1, u_2, u_3)+[r(u_1), r(u_2), r(u_3)]\\
\notag&&+(u_1, u_2)\trr' r(u_3)+c.p.+ u_1\ppr'(r(u_2), r(u_3))+c.p.\\
\notag&&-r\Big(\{u_1, u_2, u_3\}'+(u_1, u_2)\ppl' r(u_3)+c.p.\\
&&\qquad+ u_1\trl'(r(u_2), r(u_3))+c.p.\Big),
\end{eqnarray}
where c.p. means cyclic permutations with respect to elements $u_1, u_2$ and $u_3$.
\end{definition}

We denote by $\mathfrak{T}(\frkg,V)$ the set of all extending structures $\Omega(\g, V)$.
It is easy to see that  $\equiv$  and $\approx$ are equivalence relations on the set  $\mathfrak{T}(\frkg,V)$.
By the above constructions and by Theorem  \ref{thm:unifiedprod}, Theorem \ref{pt1} and Lemma \ref{mainlemma}  we obtain the following result.

\begin{thm}\label{th1}
Let $\frkg$ be a 3-Lie algebra, $E$ a vector space that contains $\frkg$ as a subspace and $V$ a complement
of $\frkg$ in $E$. Then, we get:\\
(1)Denote $\mathcal{E}\mathcal{H}^2(V,\frkg):=\mathfrak{T}(\frkg,V)/\equiv$. Then, the map
\begin{eqnarray}
\mathcal{E}\mathcal{H}^2(V,\frkg)\rightarrow Extd(E,\frkg),~~~~\overline{\Omega(\frkg,V)}\rightarrow \frkg\natural V
\end{eqnarray}
is bijective, where $\overline{\Omega(\frkg,V)}$ is the equivalence class of $\Omega(\frkg,V)$ under $\equiv$.\\
(2) Denote $\mathcal{U}\mathcal{H}^2(V,\frkg):=\mathfrak{T}(\frkg,V)/\approx$. Then, the map
\begin{eqnarray}
\mathcal{U}\mathcal{H}^2(V,\frkg)\rightarrow Extd'(E,\frkg),~~~~\overline{\overline{\Omega(\frkg,V)}}\rightarrow \frkg\natural V
\end{eqnarray}
is bijective, where $\overline{\overline{\Omega(\frkg,V)}}$ is the equivalence class of $\Omega(\frkg,V)$ under $\approx$.\\
\end{thm}

\section{Special cases of  unified products}\label{sec:specialcases}

In this section, we show some special cases of unified products and extending structures.

\subsection{Crossed products and non-abelian extension problem}
Now we give a first special case of the
unified product, namely the crossed product of 3-Lie algebras
which is related to the study of the non-abelian extension problem.

Let $\mathfrak{g}$ and $\mathfrak{h}$ be two given 3-Lie
algebras. The extension problem asks for the classification of all
extensions of $\mathfrak{h}$ by $\mathfrak{g}$, i.e. of all
3-Lie algebras ${E}$ that fit into an exact sequence
\begin{eqnarray} \label{extencros0}
\xymatrix{ 0 \ar[r] & \mathfrak{g} \ar[r]^{i} & {E}
\ar[r]^{\pi} & \mathfrak{h} \ar[r] & 0 }.
\end{eqnarray}
The classification is up to an isomorphism of 3-Lie algebras
that stabilizes $\mathfrak{g}$ and co-stabilizes $\mathfrak{h}$
and we denote by ${\mathcal E} {\mathcal P} (\mathfrak{h}, \,
\mathfrak{g})$ the isomorphism classes of all non-abelian extensions of
$\mathfrak{h}$ by $\mathfrak{g}$ up to this equivalence relation.

If $\mathfrak{g}$ is abelian, then ${\mathcal E} {\mathcal P}
(\mathfrak{h}, \, \mathfrak{g}) \cong {\mathcal{H}}^2 (\mathfrak{h}, \,
\mathfrak{g})$, where ${\mathcal{H}}^2 (\mathfrak{h}, \, \mathfrak{g})$
is the the second cohomology group in Section \ref{sec:2}. We are going to study the non-abelian case.

\begin{definition}\label{def:crossed}
Let $(\frkg,[\cdot,\cdot,\cdot]_\g)$ and $(\mathfrak{h},[\cdot,\cdot,\cdot]_{\mathfrak{h}})$ be two 3-Lie algebras.
Then $(\frkg, \mathfrak{h})$  is called a crossed system if there exits three maps
\begin{eqnarray*}
\triangleright:\mathfrak{h}\times \mathfrak{h}\times \g\to \g,\quad \ppr:\mathfrak{h}\times \g\times \g\to \g,\quad \omega:\mathfrak{h}\times \mathfrak{h}\times \mathfrak{h}\to \g,
\end{eqnarray*}
such that the direct sum space $\frkg\oplus \mathfrak{h}$ form a 3-Lie algebra under the following  bracket
\begin{eqnarray*}
&&[x_1 + u_1, x_2 + u_2, x_3 + u_3]\\
&=&[x_1, x_2, x_3] +\omega(u_1, u_2, u_3)+(u_1, u_2)\trr x_3 +(u_2, u_3)\trr x_1+ (u_3, u_1)\trr x_2\\
&&+u_3\ppr (x_1, x_2)+ u_1\ppr(x_2, x_3) +u_2\ppr (x_3, x_1)+[u_1, u_2, u_3].
\end{eqnarray*}
where $x_i\in \frkg$ and $u_i\in \mathfrak{h}$,  the bracket in $\g$ and $\mathfrak{h}$ is written as $[x_1, x_2, x_3]$ and $[u_1, u_2, u_3]$  for simplicity. This 3-Lie algebra is called crossed product of $\g$ and $\mathfrak{h}$ which will be  denoted   by $\frkg\#^\omega_{\trr,\ppr} \mathfrak{h}$.
\end{definition}

From the above  Definition \ref{def:crossed}, it is easy to see that a crossed system is a special case of  extending datum  when the maps
 $\triangleleft:V\times \g\times \g\to V$ and $\ppl:V\times V\times \g \to V$ are zero. Thus by letting these two maps to be zero in the conditions of Theorem  \ref{thm:unifiedprod},  we obtain the following conditions in Theorem \ref{thm:crossed} of  a crossed system to be  a crossed product.

\begin{thm}\label{thm:crossed}
Let $(\frkg,[\cdot,\cdot,\cdot]_\g)$ and $(\mathfrak{h},[\cdot,\cdot,\cdot]_{\mathfrak{h}})$  be two 3-Lie algebras.
Then $(\frkg, \mathfrak{h})$  is  a crossed system  if and only if the following conditions hold:
\begin{eqnarray}
 v_3\ppr[(x_1, x_2) ,(y_1, y_2)]
&=&[x_1, x_2,v_3\ppr (y_1, y_2)]-[v_3\ppr (x_1, x_2), y_1, y_2],
\end{eqnarray}
\begin{eqnarray}
\notag u_1\ppr (x_2, [y_1, y_2 ,y_3])&=&[u_1\ppr (x_2, y_1),y_2, y_3] +[u_1\ppr (x_2, y_2),y_3, y_1] \\
&&+[u_1\ppr (x_2, y_3),y_1, y_2],
\end{eqnarray}
\begin{eqnarray}
\notag[x_1, x_2,(v_1,v_2)\triangleright y_3]&=&v_2\ppr\big(v_1\ppr (x_1, x_2), y_3\big)+v_1\ppr\big(v_2\ppr (x_1, x_2), y_3\big)\\
&&+(v_1, v_2)\trr[x_1, x_2, y_3],
\end{eqnarray}
\begin{eqnarray}
\notag&&u_1\ppr\big(x_2,v_1\ppr(y_2,y_3)\big)+[(u_1, v_1)\triangleright x_2, y_2 ,y_3]\\
&=&v_1\ppr\big(u_1\ppr(x_2,y_2),y_3\big)+v_1\ppr\big(y_2,u_1\ppr(x_2,y_3)\big),
\end{eqnarray}
\begin{eqnarray}
\notag(u_1, u_2)\triangleright[y_1, y_2 ,y_3]&=&[(u_1, u_2)\triangleright y_1, y_2 ,y_3]+[y_1, (u_1, u_2)\triangleright y_2 ,y_3]\\
&&+[y_1, y_2,(u_1, u_2)\triangleright y_3],
\end{eqnarray}
\begin{eqnarray}
\notag&&[v_1, v_2 ,v_3]\ppr(x_1, x_2)+[x_1, x_2,\omega(v_1, v_2 ,v_3)]\\
\notag&=&(v_2, v_3)\trr (v_1\ppr(x_1, x_2))+(v_3, v_1)\trr (v_2\ppr(x_1, x_2))+(v_1, v_2)\trr (v_3\ppr(x_1, x_2)),
\end{eqnarray}
\begin{eqnarray}
\notag u_2\ppr\big((v_1, v_2)\trr y_3, x_1\big)&=&v_2\ppr\big(y_3,(u_2 ,v_1)\trr x_1\big)+v_1\ppr\big((u_2, v_2)\trr x_1, y_3\big)\\
&&+(v_1, v_2)\trr\big(u_2\ppr(y_3,x_1)\big),
\end{eqnarray}
\begin{eqnarray}
\notag(u_1, u_2)\trr \big(v_1\ppr(y_2 ,y_3)\big)&=&[\omega(u_1, u_2 ,v_1), y_2, y_3]+[u_1, u_2 ,v_1]\ppr(y_2, y_3)\\
\notag&&+v_1\ppr\big((u_1, u_2)\trr y_2,  y_3\big)+v_1\ppr\big(y_2,(u_1, u_2)\trr y_3\big),
\end{eqnarray}
\begin{eqnarray}
\notag[(u_1, u_2),(v_1, v_2)]\trr y_3&=&(u_1, u_2)\trr((v_1, v_2)\trr y_3)-(v_1, v_2)\trr((u_1, u_2)\trr y_3)\\
&&+v_2\ppr(\omega(u_1, u_2 ,v_1), y_3)-v_1\ppr(\omega(u_1, u_2 ,v_2), y_3),
\end{eqnarray}
\begin{eqnarray}
\notag&&(u_2, [v_1, v_2 ,v_3])\trr x_1+u_2\ppr(\omega(v_1, v_2 ,v_3),x_1)\\
\notag&=&(v_2, v_3)\trr((u_2, v_1)\trr x_1)+ (v_3, v_1)\trr((u_2, v_2)\trr x_1)+(v_1, v_2)\trr((u_2, v_3)\trr x_1),
\end{eqnarray}
\begin{eqnarray}
\notag&&\omega(u_1, u_2,[v_1, v_2 ,v_3])-\omega([u_1, u_2,v_1] v_2 ,v_3)-\omega(v_2, [u_1, u_2,v_2] ,v_3)\\
\notag&&-\omega(v_1, v_2,[u_1, u_2 ,v_3])+(u_1, u_2)\trr\omega(v_1, v_2 ,v_3)-(v_2, v_3)\trr\omega(u_1, u_2,v_1)\\
&&-(v_3, v_1)\trr\omega(u_1, u_2,v_2)-(v_1, v_2)\trr\omega(u_1, u_2 ,v_3)=0.
\end{eqnarray}
\end{thm}

Now we explain how to classify the extension problem by using crossed product.  Let $\mathfrak{g}$ and $\mathfrak{h}$
be two 3-Lie algebras and we denote by ${\mathcal C}{\mathcal S}
\, (\mathfrak{h}, \, \mathfrak{g} )$ the set of all triples
$(\triangleright, \, \ppr, \, \omega)$ such that
$(\mathfrak{g}, \, \mathfrak{h}, \, \triangleright, \,
\ppr, \, \omega)$ is a crossed system of 3-Lie algebras.
First we remark that, if $(\mathfrak{g}, \, \mathfrak{h}, \,
\triangleright, \, \ppr, \, \omega)$ is a crossed system,
then the crossed product $\mathfrak{g} \#_{\triangleright,
\ppr}^\omega \, \mathfrak{h}$ is an extension of
$\mathfrak{h}$ by $\mathfrak{g}$ via
\begin{eqnarray} \label{extencros001}
\xymatrix{ 0 \ar[r] & \mathfrak{g} \ar[r]^{i_{}} &
 \mathfrak{g}
\#_{\triangleright, \ppr}^\omega \, \mathfrak{h}
\ar[r]^{\pi_{}} & \mathfrak{h} \ar[r] & 0 }
\end{eqnarray}
where $i_{} (x) = (x, 0)$ and $\pi_{} (x,
u) = u$ are the canonical maps. Conversely, any extension $E$ of $\mathfrak{h}$ by
$\mathfrak{g}$ is equivalent to a crossed product extension of the
form \eqref{extencros001}. Thus, the classification of all
extensions of $\mathfrak{h}$ by $\mathfrak{g}$ reduces to the
classification of all crossed products $\mathfrak{g}
\#_{\triangleright, \ppr}^\omega \, \mathfrak{h}$ associated
to all crossed systems of 3-Lie algebras $(\mathfrak{g},
\mathfrak{h}, \triangleright, \ppr, \omega)$.

By lemma \ref{mainlemma}, in the special case of crossed systems, we obtain the
following

\begin{definition}
Two crossed systems  $(\triangleright,\,
\ppr, \, \omega)$ and $(\triangleright ',\, \ppr ',
\, \omega')$ of ${\mathcal C}{\mathcal S} \, (\mathfrak{h}, \,
\mathfrak{g} )$ are called \emph{cohomologous} and we denote this by
$(\triangleright,\, \ppr, \, \omega) \approx (\triangleright
',\, \ppr ', \, \omega')$ if there exists a linear map $r:
\mathfrak{h} \to \mathfrak{g}$ satisfying the following  conditions for all $x_i \in \g, u_i\in \mathfrak{h}$:
\begin{eqnarray}
u_1\ppr (x_2, x_3)&=&u_1\ppr' (x_2, x_3)+[x_2, x_3, r(u_1)],\\
(u_1, u_2)\triangleright x_3&=&(u_1,u_2)\trr' x_3+[r(u_1), r(u_2), x_3],\\
\notag\omega(u_1, u_2, u_3)&=&\omega'(u_1,u_2, u_3)+[r(u_1), r(u_2), r(u_3)]-r([u_1, u_2, u_3])\\
&&+(u_1,u_2)\trr' r(u_3)+c.p.+\nu(u_1)\ppr'(r(u_2), r(u_3))+c.p.
\end{eqnarray}
\end{definition}

Note that $(\triangleright,\,\ppr, \, \omega) \approx (\triangleright ',\, \ppr', \, \omega')$ if and only if there exists $\psi : \mathfrak{g}\#_{\triangleright, \ppr}^\omega \, \mathfrak{h} \to\mathfrak{g} \#_{\triangleright ', \ppr '}^{f'} \,\mathfrak{h}$ an isomorphism of 3-Lie algebras that stabilizes $\mathfrak{g}$ and co-stabilizes $\mathfrak{h}$.
Thus we obtain the classifying result to the extension problem in the non-abelian case:

\begin{pro}
Let $\mathfrak{g}$ and $\mathfrak{h}$ be two 3-Lie
algebras. Then  $\approx$ is an equivalence relation on the set
${\mathcal C}{\mathcal S} \, (\mathfrak{h}, \, \mathfrak{g} )$ of
all crossed systems and the map
$$
{\mathcal{NH}}^2 (\mathfrak{h}, \, \mathfrak{g}) :=
{\mathcal C}{\mathcal S} \, (\mathfrak{h}, \, \mathfrak{g} )/
\approx \,\, \longrightarrow {\mathcal E} {\mathcal P}
(\mathfrak{h}, \, \mathfrak{g}), \qquad
\overline{(\triangleright,\, \ppr, \, \omega)} \mapsto
\mathfrak{g} \#_{\triangleright, \ppr}^\omega \, \mathfrak{h}
$$
is a bijection,  where $\overline{(\triangleright,\, \ppr, \, \omega)}$  is the equivalence class of $(\triangleright,\, \ppr, \, \omega)$ under $\approx$.
\end{pro}

\begin{ex}\cite{SMT}
Let $\frkg$ be the three dimensional 3-Lie algebra defined with respect to a basis
$\{x_1, x_2, x_3\}$ by the skew-symmetric bracket $[x_1, x_2, x_3] = x_1$,
and let $\frkh$  be the same 3-Lie algebra
which we consider with respect to basis $\{v_1, v_2, v_3\}$, that is $[v_1, v_2, v_3] = v_1$.
Then we have a crossed system given by:
\begin{eqnarray*}
(u_i,u_j)\trr x_k  &=& 0,\\
u_1\ppr(x_1, x_2) &=& 0,\quad u_1\ppr(x_1, x_3) = 0,\quad u_1\ppr(x_2, x_3) = r_1x_1,\\
u_2\ppr(x_1, x_2) &=& 0,\quad u_2\ppr(x_1, x_3) = 0,\quad u_2\ppr(x_2, x_3) = r_2x_1,\\
u_3\ppr(x_1, x_2) &=& 0,\quad u_3\ppr(x_1, x_3) = 0,\quad u_3\ppr(x_2, x_3) = r_3x_1,\\
\omega(u_1, u_2, u_3) &=& -r_1x_1,
\end{eqnarray*}
where $r_i$ are arbitrary parameters.
\end{ex}

\subsection{Matched pair of 3-Lie algebras}

The concept of a matched pair of Lie algebras
was introduced in quantum group theory \cite[Theorem 4.1]{majid} and  in Poisson Lie group thory \cite[Theorem 3.9]{LW}.  Now we
introduce the concept of a matched pair of 3-Lie algebras.

\begin{definition}\label{def:matchedpair}
Let $(\frkg,[\cdot,\cdot,\cdot]_\g)$ and $(\mathfrak{h},[\cdot,\cdot,\cdot]_{\mathfrak{h}})$  be two 3-Lie algebras.
Then $(\frkg, \mathfrak{h})$  is called a matched pair if there exits four linear maps
\begin{eqnarray*}
\triangleright:\mathfrak{h}\times \mathfrak{h}\times \g\to \g,\quad \triangleleft:\mathfrak{h}\times \g\times \g\to \mathfrak{h}, \\
\ppr:\mathfrak{h}\times \g\times \g\to \g,\quad \ppl:\mathfrak{h}\times \mathfrak{h}\times \g \to \mathfrak{h},
\end{eqnarray*}
such that the direct sum space $\frkg\oplus \mathfrak{h}$ form a 3-Lie algebra under the following bracket:
\begin{eqnarray*}
&&[x_1 + u_1, x_2 + u_2, x_3 + u_3]\\
&=&[x_1, x_2, x_3]+ (u_1, u_2)\trr x_3 +(u_2, u_3)\trr x_1+ (u_3, u_1)\trr x_2\\
&&+u_3\ppr (x_1, x_2)+ u_1\ppr(x_2, x_3) +u_2\ppr (x_3, x_1)\\
&&+[u_1, u_2, u_3]+ u_1\trl(x_2, x_3)+ u_2\trl(x_3, x_1) + u_3 \trl(x_1, x_2) \\
&&+ (u_2, u_3)\ppl x_1+ (u_3, u_1)\ppl x_2+ (u_1, u_2)\ppl x_3,
\end{eqnarray*}
where $x_1, x_2, x_3 \in \frkg$ and $u_1, u_2, u_3 \in \mathfrak{h}$.
This 3-Lie algebra is called the bicrossed product of $\frkg$ and $\mathfrak{h}$ . We will denoted it  by $\frkg\bowtie \mathfrak{h}$.
\end{definition}

From the above  Definition \ref{def:matchedpair}, one see that a bicrossed product  is a special case of  a unified product  when the map
 $\omega$ is zero.
 Thus from the conditions in Theorem  \ref{thm:unifiedprod},  we obtain the following Theorem \ref{thm:matchedpair} for a matched pair.

\begin{thm}\label{thm:matchedpair}
Let $(\frkg,[\cdot,\cdot,\cdot])$ and $(\mathfrak{h},[\cdot,\cdot,\cdot])$  be two 3-Lie algebras,
Then $(\frkg, \mathfrak{h})$  is  a matched pair if and only if the following compatibility conditions hold:
\begin{eqnarray}\label{eq:matched01}
v_3\trl[(x_1, x_2),(y_1, y_2)]&=&(v_3\trl(x_1, x_2))\trl(y_1, y_2)-(v_3\trl(y_1, y_2))\trl(x_1, x_2),
\end{eqnarray}
\begin{eqnarray}\label{eq:matched02}
\notag u_1\triangleleft (x_2, [y_1, y_2 ,y_3])&=&(u_1\triangleleft (x_2, y_1))\trl (y_2, y_3)+(u_1\triangleleft (x_2, y_2))\trl (y_3, y_1)\\
&&+(u_1\triangleleft (x_2, y_3))\trl (y_1, y_2),
\end{eqnarray}
\begin{eqnarray}
\notag v_3\ppr[(x_1, x_2) ,(y_1, y_2)]&=&[x_1, x_2,v_3\ppr (y_1, y_2)]+(v_3\trl(y_1, y_2)) \ppr(x_1, x_2)\\
&&-[v_3\ppr (x_1, x_2), y_1, y_2]-(v_3\trl(x_1, x_2)) \ppr(y_1, y_2),
\end{eqnarray}
\begin{eqnarray}
\notag  u_1\ppr (x_2, [y_1, y_2 ,y_3])&=&[u_1\ppr (x_2, y_1),y_2, y_3] +(u_1\triangleleft (x_2, y_1))\ppr (y_2, y_3)\\
\notag&&+[u_1\ppr (x_2, y_2),y_3, y_1] +(u_1\triangleleft (x_2, y_2))\ppr (y_3, y_1)\\
&&+[u_1\ppr (x_2, y_3),y_1, y_2] +(u_1\triangleleft (x_2, y_3))\ppr (y_1, y_2),
\end{eqnarray}
\begin{eqnarray}
\notag&&[x_1, x_2,(v_1,v_2)\triangleright y_3]+((v_1,v_2)\ppl y_3)\ppr (x_1, x_2)\\
\notag&=&v_2\ppr\big(v_1\ppr (x_1, x_2), y_3\big)+\big(v_1\trl (x_1, x_2), v_2\big)\trr y_3\\
&&+v_1\ppr \big(v_2\ppr (x_1, x_2), y_3\big)+\big(v_1, v_2\trl(x_1, x_2)\big)\trr y_3+(v_1, v_2)\trr[x_1, x_2, y_3],
\end{eqnarray}
\begin{eqnarray}
\notag((v_1,v_2)\ppl y_3)\trl (x_1, x_2)&=&v_2\trl\big(v_1\ppr (x_1, x_2), y_3\big)+\big(v_1\trl (x_1, x_2), v_2\big)\ppl y_3\\
\notag&&+v_1\trl\big(v_2\ppr (x_1, x_2), y_3\big)+\big(v_1, v_2\trl(x_1, x_2)\big)\ppl y_3\\
&&+(v_1, v_2)\ppl[x_1, x_2, y_3],
\end{eqnarray}
\begin{eqnarray}
\notag&&u_1\ppr\big(x_2,v_1\ppr(y_2,y_3)\big)+\big(v_1,u_1\trl(x_2,y_3)\big)\trr y_2\\
\notag&&+[(u_1,v_1)\trr x_2,y_2,y_3]+((u_1,v_1)\ppl x_2)\ppr(y_2,y_3)\\
\notag&=&v_1\ppr\big(u_1\ppr(x_2,y_2),y_3\big)+\big(v_1, u_1\trl(x_2,y_2)\big)\trr y_3\\
&&+v_1\ppr\big(y_2,u_1\ppr(x_2,y_3)\big)+\big(u_1,v_1\trl(y_2,y_3)\big)\trr x_2,
\end{eqnarray}
\begin{eqnarray}
\notag&&u_1\trl\big(x_2,v_1\ppr(y_2,y_3)\big)+((u_1,v_1)\ppl x_2)\trl(y_2,y_3)+\big(v_1,u_1\trl(x_2,y_3)\big)\ppl y_2\\
\notag&=&\big(u_1,v_1\trl(y_2,y_3)\big)\ppl x_2+v_1\trl\big(u_1\ppr(x_2,y_2),y_3\big)\\
&&+\big(v_1, u_1\trl(x_2,y_2)\big)\ppl y_3+v_1\trl\big(y_2,u_1\ppr(x_2,y_3)\big),
\end{eqnarray}
\begin{eqnarray}
\notag(u_1, u_2)\triangleright[y_1, y_2 ,y_3]&=&[(u_1, u_2)\triangleright y_1, y_2 ,y_3]+((u_1, u_2)\ppl y_1)\ppr(y_2 ,y_3)\\
\notag&&+[y_1, (u_1, u_2)\triangleright y_2 ,y_3]+((u_1, u_2)\ppl y_2)\ppr(y_3 ,y_1)\\
&&+[y_1, y_2,(u_1, u_2)\triangleright y_3]+((u_1, u_2)\ppl y_3)\ppr(y_1 ,y_2),
\end{eqnarray}
\begin{eqnarray}
\notag(u_1, u_2)\ppl[y_1, y_2 ,y_3]&=&((u_1, u_2)\ppl y_1)\trl(y_2 ,y_3)+((u_1, u_2)\ppl y_2)\trl(y_3 ,y_1)\\
&&+((u_1, u_2)\ppl y_3)\trl(y_1 ,y_2),
\end{eqnarray}
\begin{eqnarray}
\notag[v_1, v_2 ,v_3]\ppr(x_1, x_2)&=&(v_2, v_3)\trr (v_1\ppr(x_1, x_2))+(v_3, v_1)\trr (v_2\ppr(x_1, x_2))\\
&&+(v_1, v_2)\trr (v_3\ppr(x_1, x_2)),
\end{eqnarray}
\begin{eqnarray}
\notag[v_1, v_2 ,v_3]\trl(x_1, x_2)&=&(v_2, v_3)\ppl (v_1\ppr(x_1, x_2))+[v_1\trl(x_1, x_2), v_2, v_3]\\
\notag&&+(v_3, v_1)\ppl (v_2\ppr(x_1, x_2))+[v_2\trl(x_1, x_2), v_3, v_1]\\
&&+(v_1, v_2)\ppl (v_3\ppr(x_1, x_2))+[v_3\trl(x_1, x_2), v_1, v_2].
\end{eqnarray}
\begin{eqnarray}
\notag&&u_2\ppr\big((v_1, v_2)\trr y_3, x_1\big)+\big(u_2,(v_1, v_2)\ppl y_3\big)\trr x_1\\
\notag&=&v_2\ppr\big(y_3,(u_2 ,v_1)\trr x_1\big)+\big((u_2 ,v_1)\ppl x_1, v_2\big)\trr y_3+v_1\ppr\big((u_2, v_2)\trr x_1, y_3\big)\\
&&+\big(v_1, (u_2, v_2)\ppl x_1)\trr y_3+(v_1, v_2)\trr\big(u_2\ppr(y_3,x_1)\big),
\end{eqnarray}
\begin{eqnarray}
\notag&&u_2\trl\big((v_1, v_2)\trr y_3, x_1\big)+\big(u_2,(v_1, v_2)\ppl y_3\big)\ppl x_1\\
\notag&=&v_2\trl\big(y_3,(u_2 ,v_1)\trr x_1\big)+\big((u_2 ,v_1)\ppl x_1, v_2\big)\ppl y_3\\
\notag&&+v_1\trl\big((u_2, v_2)\trr x_1, y_3\big)+\big(v_1, (u_2, v_2)\ppl x_1)\ppl y_3\\
&&+(v_1, v_2)\ppl \big(u_2\ppr(y_3,x_1)\big)+[v_1, v_2,u_3\trl(y_3,x_1)].
\end{eqnarray}
\begin{eqnarray}
\notag(u_1, u_2)\trr \big(v_1\ppr(y_2 ,y_3)\big)&=&[u_1, u_2 ,v_1]\ppr(y_2, y_3)+v_1\ppr\big((u_1, u_2)\trr y_2,  y_3\big)\\
\notag&&+\big(v_1, (u_1, u_2)\ppl y_2\big)\trr y_3+v_1\ppr\big(y_2,(u_1, u_2)\trr y_3\big)\\
&&+\big((u_1, u_2)\ppl y_3,v_1\big)\trr y_2,
\end{eqnarray}
\begin{eqnarray}
\notag&&(u_1, u_2)\ppl \big(v_1\ppr(y_2 ,y_3)\big)+[u_1, u_2, v_1\trl(y_2 ,y_3)]\\
\notag&=&[u_1, u_2 ,v_1]\trl(y_2, y_3)+v_1\trl\big((u_1, u_2)\trr y_2,  y_3\big)+\big(v_1, (u_1, u_2)\ppl y_2\big)\ppl y_3\\
&&+v_1\trl\big(y_2,(u_1, u_2)\trr y_3\big)+\big((u_1, u_2)\ppl y_3,v_1\big)\ppl y_2,
\end{eqnarray}
\begin{eqnarray}
\notag[(u_1, u_2), (v_1, v_2)]\ppl y_3&=&(u_1, u_2)\ppl((v_1, v_2)\trr y_3)+[u_1, u_2,(v_1, v_2)\ppl y_3]\\
&&-(v_1, v_2)\ppl((u_1, u_2)\trr y_3)-[v_1, v_2,(u_1, u_2)\ppl y_3],
\end{eqnarray}
\begin{eqnarray}\label{eq:matched18}
[(u_1, u_2),(v_1, v_2)]\trr y_3&=&(u_1, u_2)\trr((v_1, v_2)\trr y_3)-(v_1, v_2)\trr((u_1, u_2)\trr y_3),
\end{eqnarray}
\begin{eqnarray}\label{eq:matched19}
\notag(u_2, [v_1, v_2 ,v_3])\trr x_1&=&(v_2, v_3)\trr((u_2, v_1)\trr x_1)+ (v_3, v_1)\trr((u_2, v_2)\trr x_1)\\
&&+(v_1, v_2)\trr((u_2, v_3)\trr x_1),
\end{eqnarray}
\begin{eqnarray}
\notag (u_2, [v_1, v_2 ,v_3])\ppl x_1&=&(v_2, v_3)\trr((u_2, v_1)\ppl x_1) +[(u_2, v_1)\ppl x_1,v_2, v_3]\\
\notag&&+ (v_3, v_1)\trr((u_2, v_2)\ppl x_1) +[(u_2, v_2)\ppl x_1,v_3, v_1]\\
&&+ (v_1, v_2)\trr((u_2, v_3)\ppl x_1) +[(u_2, v_3)\ppl x_1,v_1, v_2].
\end{eqnarray}
The equations \eqref{eq:matched01}-\eqref{eq:matched02} means that $(\mathfrak{h},\trl)$ is a representation of $\g$ and the equations\eqref{eq:matched18}-\eqref{eq:matched19} means that $(\g,\trr)$ is a representation of $\mathfrak{h}$.
\end{thm}

\begin{ex}
Let $\mathfrak{g}$ be the  three dimensional 3-Lie algebra  with non-zero bracket under the basis $\{x_1, x_2, x_3\}$,
$$\left[x_{1}, x_{2}, x_{3}\right]=x_{3},$$
 and  $\mathfrak{h}$ be a  five dimension 3-Lie algebra with basis $\{u_1, u_2, u_3, u_4, u_5\}$.
 $$\left[u_{2}, u_{3}, u_{4}\right]=u_{1}.$$
  We define the following extending datum:
 $$
\begin{gathered}
u_5\ppr (x_{1}, x_{2})=  x_{3}, \\
(u_{3}, u_{5})\ppl x_{2}=  u_{1},\quad (u_{3}, u_{5})\ppl x_{1}= - u_{2},\\
(u_{4},u_{5})\ppl x_{2}= - u_{1}, (u_{4}, u_{5})\ppl x_{1}= u_{2},\\
(u_{1}, u_{2})\trr  x_{1}=-x_{3},\quad (u_{1}, u_{3})\trr x_{1}=x_{2},\quad (u_{2}, u_{3})\trr x_{1}=-x_{1}, \\
u_{3}\trl (x_{1}, x_{2})=-  u_{3}-  u_{4},\quad u_{3}\trl (x_{2}, x_{3})=-  u_{1},\quad u_{3}\trl( x_{1}, x_{3})=  u_{2}, \\
u_{4}\trl (x_{1}, x_{2})=-  u_{3}+  u_{4},\quad u_{5}\trl (x_{1}, x_{2})=-  u_{5}, \\
\end{gathered}
$$
Then one checks that this extending datum satisfying the matched pair conditions in the above  Theorem \ref{thm:matchedpair}.
Therefore  we obtain an eight dimensional 3-Lie algebra $E=\mathfrak{g}\bowtie \mathfrak{h}$ defined with respect to a basis
$\{x_1, x_2, x_3,u_1, u_2, u_3,  u_4, u_5\}$ by the skew-symmetric bracket
$$
\begin{gathered}
\left[x_{1}, x_{2}, x_{3}\right]=x_{3},\quad \left[u_{2}, u_{3}, u_{4}\right]=u_{1},\\
\left[x_{1}, x_{2}, u_{5}\right]=  x_{3}-  u_{5}, \\
\left[u_{3}, x_{2}, u_{5}\right]=-  u_{1},\quad \left[u_{3}, x_{1}, u_{5}\right]=  u_{2},\\
\left[u_{4}, x_{2}, u_{5}\right]=  u_{1},\quad\left[u_{4}, x_{1}, u_{5}\right]=-  u_{2},\\
\left[u_{1}, u_{2}, x_{1}\right]=-x_{3},\quad\left[u_{1}, u_{3}, x_{1}\right]=x_{2},\quad\left[u_{2}, u_{3}, x_{1}\right]=-x_{1}, \\
\left[u_{3}, x_{1}, x_{2}\right]=-  u_{3}-  u_{4},\quad \left[u_{3}, x_{2}, x_{3}\right]=-  u_{1},\quad\left[u_{3}, x_{1}, x_{3}\right]=  u_{2}, \\
\left[u_{4}, x_{1}, x_{2}\right]=-  u_{3}+  u_{4}.
\end{gathered}
$$
\end{ex}

The bicrossed product of two 3-Lie algebras is related to the so called \emph{factorization problem}, which  can be stated as follows:
\emph{Let $\mathfrak{g}$ and $\mathfrak{h}$ be two given 3-Lie
algebras. Describe and classify all 3-Lie algebras ${E}$ that
factorize through $\mathfrak{g}$ and $\mathfrak{h}$, i.e. ${E}$
contains $\mathfrak{g}$ and $\mathfrak{h}$ as 3-Lie subalgebras
such that ${E} = \mathfrak{g} + \mathfrak{h}$ and $\mathfrak{g}
\cap \mathfrak{h} = [0]$.} We use Theorem \ref{pt1} to prove the following:

\begin{pro}\label{bicrfactor}
A 3-Lie algebra ${E}$ factorizes through $\mathfrak{g}$ and
$\mathfrak{h}$ if and only if there exists a matched pair of
3-Lie algebras $(\mathfrak{g}, \mathfrak{h}, \, \triangleleft,
\, \triangleright, \, \leftharpoonup, \, \rightharpoonup)$ such
that $ {E} \cong \mathfrak{g} \bowtie \mathfrak{h}$.
\end{pro}

\begin{proof}
We know that any bicrossed product $\mathfrak{g}
\bowtie \mathfrak{h}$ factorizes through $\mathfrak{g} \cong
\mathfrak{g} \times [0]$ and $\mathfrak{h} \cong [0]\times
\mathfrak{h}$. Conversely, assume that ${E}$ factorizes through
$\mathfrak{g}$ and $\mathfrak{h}$. Let $p: {E} \to \mathfrak{g}$
be the natural projection of ${E}$ on $\mathfrak{g}$, i.e. $p
(x + u) := x$, for all $x\in \mathfrak{g}$ and $u \in
\mathfrak{h}$. Now, we apply Theorem \ref{pt1} for $V: = {\rm
Ker}(p) = \mathfrak{h}$. Since $V$ is a 3-Lie subalgebra of ${E}$, the map $\omega$  is the trivial map and the extending structure $\Omega(\mathfrak{g},V) = \bigl(\triangleleft, \, \triangleright, \, \leftharpoonup, \,
\rightharpoonup, \, [-, \, -]\bigl)$ constructed in the
proof of Theorem \ref{pt1} is precisely a matched pair of 3-Lie
algebra. Thus the unified product $\mathfrak{g}\natural  V $ is
the bicrossed product $\mathfrak{g} \bowtie \mathfrak{h}$. Explicitly, the matched pair $(\mathfrak{g},
\mathfrak{h}, \, \triangleleft, \, \triangleright \, \leftharpoonup, \, \rightharpoonup)$ is
given by:
\begin{eqnarray*}
\triangleright:\mathfrak{h}\times \mathfrak{h}\times \g\to \g,&&\quad (u_1, u_2)\triangleright x_3=p([u_1, u_2, x_3]), \\
\ppr:\mathfrak{h}\times \g\times \g\to \g,&&\quad u_1\ppr (x_2, x_3)=p([u_1, x_2, x_3]),\\
 \triangleleft:\mathfrak{h}\times \g\times \g\to \mathfrak{h}, &&\quad u_1\triangleleft (x_2, x_3)=[u_1, x_2, x_3]-p([u_1, x_2, x_3]),\\
 \ppl:\mathfrak{h}\times \mathfrak{h}\times \g \to \mathfrak{h},&&\quad (u_1, u_2)\ppl x_3=[u_1, u_2, x_3]-p([u_1, u_2, x_3]),
\end{eqnarray*}
for all $u_i \in \mathfrak{h}$ and $x_i \in \mathfrak{g}$.
\end{proof}


\subsection{Classifying complements for 3-Lie algebras} \label{complements}

This subsection is devoted to the classifying complements problem. Let
$\mathfrak{g} \subseteq {E}$ be a 3-Lie subalgebra of ${E}$. A
3-Lie subalgebra $\mathfrak{h}$ of ${E}$ is called a
\emph{complement} of $\mathfrak{g}$ in ${E}$ (or a
\emph{$\mathfrak{g}$-complement} of ${E}$) if ${E} = \mathfrak{g}
+ \mathfrak{h}$ and $\mathfrak{g} \cap \mathfrak{h} = [0]$. If
$\mathfrak{h}$ is a complement of $\mathfrak{g}$ in ${E}$,
then we have  $ {E} \cong \mathfrak{g} \bowtie
\mathfrak{h}$, where $\mathfrak{g} \bowtie \mathfrak{h}$ is the
bicrossed product associated to the matched pair of the factorization ${E} = \mathfrak{g} + \mathfrak{h}$, see Proposition \ref{bicrfactor}.

We denote by ${\mathcal F} (\mathfrak{g}, \, {E})$ the (possibly
empty) isomorphism classes of all $\mathfrak{g}$-complements of
${E}$. The \emph{factorization index} of $\mathfrak{g}$ in $E$
is defined by $[E : \mathfrak{g}]:= |\, {\mathcal F}
(\mathfrak{g}, \, E) \,|$.

\begin{definition} \label{deformaplie}
Let $(\mathfrak{g}, \, \mathfrak{h}, \, \triangleright, \,
\triangleleft, \, \leftharpoonup, \, \rightharpoonup)$ be a
matched pair of 3-Lie algebras. A linear map $r: \mathfrak{h}
\to \mathfrak{g}$ is called a \emph{deformation map} of the
matched pair $(\mathfrak{g}, \mathfrak{h}, \triangleright,
\triangleleft, \, \leftharpoonup, \, \rightharpoonup)$ if the
following condition holds for any $u_i \in \mathfrak{h}$:
\begin{eqnarray}\label{factLie}
\notag&&r[u_1, u_2, u_3]-[r(u_1), r(u_2), r(u_3)]\\
\notag&=&(u_1, u_2)\trr r(u_3)+c.p.+ u_1\ppr(r(u_2), r(u_3))+c.p.\\
&&-r\Big((u_1, u_2)\ppl r(u_3)+c.p.+ u_1\trl(r(u_2), r(u_3))+c.p.\Big).
\end{eqnarray}
\end{definition}

We denote by ${\mathcal D}{\mathcal M} \, (\mathfrak{h},
\mathfrak{g} \, | \, (\triangleright, \triangleleft,
\leftharpoonup, \, \rightharpoonup) )$ the set of all deformation
maps of the matched pair $(\mathfrak{g}, \mathfrak{h},
\triangleright, \triangleleft, \leftharpoonup, \rightharpoonup)$.
The trivial map $r(x) = 0$, for all $x \in \mathfrak{h}$, is of
course a deformation map. The right hand side of \eqref{factLie}
measures how far $r: \mathfrak{h} \to \mathfrak{g}$ is from being
a 3-Lie algebra map. Using this concept of deformation map, we introduce the following
deformation of 3-Lie algebras.

\begin{thm}\label{deforLie}
Let $\mathfrak{g}$ be a 3-Lie subalgebra of ${E}$,
$\mathfrak{h}$ a given $\mathfrak{g}$-complement of ${E}$ and $r:
\mathfrak{h} \to \mathfrak{g}$ a deformation map of the associated
canonical matched pair $(\mathfrak{g}, \mathfrak{h},
\triangleright, \triangleleft, \leftharpoonup, \rightharpoonup)$.

$(1)$ Let $f_{r}: \mathfrak{h} \to {E} = \mathfrak{g} \bowtie
\mathfrak{h}$ be the linear map defined for any $u \in
\mathfrak{h}$ by:
$$f_{r}(u) = (r(u),\, u).$$
Then $\widetilde{\mathfrak{h}} : = {\rm Im}(f_{r})$ is a
$\mathfrak{g}$-complement of ${E}$.

$(2)$ Let $\mathfrak{h}_{r} := \mathfrak{h}$, as a vector space, with
the new bracket defined for any $u_i\in \mathfrak{h}$ by:
\begin{equation}\label{rLiedef}
[u_1, u_2, u_3]_{r} := [u_1, u_2, u_3] + (u_1, u_2)\ppl r(u_3)+c.p.+ u_1\trl(r(u_2), r(u_3))+c.p.
\end{equation}
Then $\mathfrak{h}_{r} $ is a 3-Lie algebra which is called the $r$-deformation of $\mathfrak{h}$.
Furthermore, $\mathfrak{h}_{r} \cong \widetilde{\mathfrak{h}}$, as
3-Lie algebras.
\end{thm}

\begin{proof}
$(1)$ First we will prove that $\widetilde{\mathfrak{h}}
= [\bigl(r(u),\, u \bigl) ~|~ u \in \mathfrak{h}]$ is a 3-Lie
subalgebra of $\mathfrak{g} \bowtie \mathfrak{h} = {E} $. Indeed,
for all $u_i\in \mathfrak{h}$ we have:
\begin{eqnarray*}
&&\bigl[ (r(u_1), u_1), (r(u_2), u_2), (r(u_3), u_3)  \bigl]\\
&{=}& \Bigl(\bigl[r(u_1),r(u_2),r(u_3)\bigl] + (u_1, u_2)\trr r(u_3)+c.p.+ u_1\ppr(r(u_2), r(u_3))+c.p.,\\
&&\qquad[u_1, \, u_2, \, u_3] +(u_1, u_2)\ppl r(u_3)+c.p.+ u_1\trl(r(u_2), r(u_3))+c.p.\Bigl)\\
&\stackrel{\eqref{factLie}}{=}& \Bigl(r \big([u_1, \, u_2, \, u_3]+(u_1, u_2)\ppl r(u_3)+c.p.+ u_1\trl(r(u_2), r(u_3))+c.p.\big),\\
&&\qquad[u_1, \, u_2, \, u_3] +(u_1, u_2)\ppl r(u_3)+c.p.+ u_1\trl(r(u_2), r(u_3))+c.p.\Bigl).
\end{eqnarray*}
Thus $\bigl[ (r(u_1), u_1), (r(u_2), u_2), (r(u_3), u_3)  \bigl]\in
\widetilde{\mathfrak{h}}$. Moreover, it is straightforward to see
that $\mathfrak{g} \, \cap \,  \widetilde{\mathfrak{h}} = [0]$
and $(x,\, u) = \bigl(x- r(u), \, 0\bigl) + \bigl(r(u),\, u
\bigl) \in \mathfrak{g} + \widetilde{\mathfrak{h}}$ for all $x \in
\mathfrak{g}$, $u \in \mathfrak{h}$. Here, we view $\mathfrak{g}
\cong \mathfrak{g} \times [0]$ as a subalgebra of $\mathfrak{g}
\bowtie \mathfrak{h}$. Therefore, $\widetilde{\mathfrak{h}}$ is a
$\mathfrak{g}$-complement of ${E} = \mathfrak{g} \bowtie
\mathfrak{h}$.

$(2)$ We denote by $\widetilde{f_{r}} : \mathfrak{h} \to
\widetilde{\mathfrak{h}}$ the linear isomorphism induced by
$f_{r}$. We will prove that $\widetilde{f_{r}}$ is also a 3-Lie
algebra map if we consider on $\mathfrak{h}$ the bracket given by
\eqref{rLiedef}. In fact, for any $u_i\in \mathfrak{h}$ we
have:
\begin{eqnarray*}
&&\widetilde{f_{r}}\bigl([u_1, \, u_2, \, u_3]_{r}\bigl)\\
&\stackrel{\eqref{rLiedef}}{=}&
\widetilde{f_{r}}\bigl([u_1, u_2, u_3] + (u_1, u_2)\ppl r(u_3)+c.p.+ u_1\trl(r(u_2), r(u_3))+c.p.\bigl)\\
&{=}& \Bigl({r \bigl([u_1, u_2, u_3] + (u_1, u_2)\ppl r(u_3)+c.p.+ u_1\trl(r(u_2), r(u_3))+c.p.\bigl)},\,\\
&&\qquad[u_1, u_2, u_3] + (u_1, u_2)\ppl r(u_3)+c.p.+ u_1\trl(r(u_2), r(u_3))+c.p.\Bigl)\\
&\stackrel{\eqref{factLie}}{=}& \Bigl(\bigl[r(u_1),r(u_2),r(u_3)\bigl] + (u_1, u_2)\trr r(u_3)+c.p.+ u_1\ppr(r(u_2), r(u_3))+c.p.,\\
&&\qquad[u_1, \, u_2, \, u_3] +(u_1, u_2)\ppl r(u_3)+c.p.+ u_1\trl(r(u_2), r(u_3))+c.p.\Bigl)\\
&{=}& [(r(u_1),\, u_1), \, (r(u_2), \, u_2), \, (r(u_3), \, u_3)]
= [\widetilde{f_{r}}(u_1), \, \widetilde{f_{r}}(u_2), \, \widetilde{f_{r}}(u_3)].
\end{eqnarray*}
Therefore, $\mathfrak{h}_{r}$ is a 3-Lie algebra and the proof is  finished.
\end{proof}

The following is the converse of Theorem \ref{deforLie}. It proves that
all $\mathfrak{g}$-complements of ${E}$ are $r$-deformations of a
given complement.

\begin{thm} \label{descrierecomlie}
Let $\mathfrak{g}$ be a 3-Lie subalgebra of ${E}$,
$\mathfrak{h}$ a given $\mathfrak{g}$-complement of ${E}$ with the
associated canonical matched pair of 3-Lie algebras
$(\mathfrak{g}, \mathfrak{h}, \triangleright, \triangleleft, ,
\leftharpoonup, \rightharpoonup)$. Then $\overline{\mathfrak{h}}$
is a $\mathfrak{g}$-complement of ${E}$ if and only if there
exists an isomorphism of 3-Lie algebras $\overline{\mathfrak{h}}
\cong \mathfrak{h}_{r}$, for some deformation map $r: \mathfrak{h}
\to \mathfrak{g}$ of the matched pair $(\mathfrak{g},
\mathfrak{h}, \triangleright, \triangleleft, \leftharpoonup,
\rightharpoonup)$.
\end{thm}

\begin{proof}
Let $\overline{\mathfrak{h}}$ be an arbitrary
$\mathfrak{g}$-complement of ${E}$. Since ${E} = \mathfrak{g}
\oplus \mathfrak{h} = \mathfrak{g} \oplus \overline{\mathfrak{h}}$
we can find four $k$-linear maps:
$$
s: \mathfrak{h} \to \mathfrak{g}, \quad \nu: \mathfrak{h} \to
\overline{\mathfrak{h}}, \quad t:\overline{\mathfrak{h}} \to
\mathfrak{g}, \quad \varpi: \overline{\mathfrak{h}} \to \mathfrak{h}
$$
such that for all $u \in \mathfrak{h}$ and $v \in
\overline{\mathfrak{h}}$ we have:
\begin{equation} \label{lie111}
u = s(u) \oplus \nu(u), \qquad y = t(v) \oplus \varpi(v).
\end{equation}
It is easy to see that $\nu: \mathfrak{h} \to
\overline{\mathfrak{h}}$ is a linear isomorphism of vector spaces.
We denote by $\tilde{\nu}: \mathfrak{h} \to \mathfrak{g} \bowtie
\mathfrak{h}$ the composition:
$$
\tilde{\nu} : \, \mathfrak{h} \, \stackrel{v} {\longrightarrow} \,
\overline{\mathfrak{h}} \, \stackrel{i}{\hookrightarrow} \, {E} \,
= \,\mathfrak{g} \bowtie \mathfrak{h}.
$$
Therefore, we have $\tilde{\nu}(u) \stackrel{\eqref{lie111}}{=}
\bigl(-s(u),\, u\bigl)$, for all $u \in \mathfrak{h}$. Then we
shall prove that $r := - s$ is a deformation map and
$\overline{\mathfrak{h}} \cong \mathfrak{h}_{r}$. Indeed,
$\overline{\mathfrak{h}} = {\rm Im} (\nu) = {\rm Im} (\tilde{\nu})$ is
a 3-Lie subalgebra of ${E} = \mathfrak{g} \bowtie \mathfrak{h}$
and we have:
\begin{eqnarray*}
&&\bigl[ (r(u_1), u_1), (r(u_2), u_2), (r(u_3), u_3)  \bigl]\\
&{=}& \Bigl(\bigl[r(u_1),r(u_2),r(u_3)\bigl] + (u_1, u_2)\trr r(u_3)+c.p.+ u_1\ppr(r(u_2), r(u_3))+c.p.,\\
&&\qquad[u_1, \, u_2, \, u_3] +(u_1, u_2)\ppl r(u_3)+c.p.+ u_1\trl(r(u_2), r(u_3))+c.p.\Bigl)\\
&=&(r(v), v)
\end{eqnarray*}
for some $v\in \mathfrak{h}$. Thus, we obtain:
\begin{eqnarray}\label{lie113}
r(v) = \bigl[r(u_1),r(u_2),r(u_3)\bigl] + (u_1, u_2)\trr r(u_3)+c.p.+ u_1\ppr(r(u_2), r(u_3))+c.p.,\\
\label{lie114}
 v = [u_1, \, u_2, \, u_3] +(u_1, u_2)\ppl r(u_3)+c.p.+ u_1\trl(r(u_2), r(u_3))+c.p.
\end{eqnarray}

By applying $r$ to \eqref{lie114} it follows
that $r$ is a deformation map of the matched pair $(\mathfrak{g},
\mathfrak{h}, \triangleright, \triangleleft, \leftharpoonup,
\rightharpoonup)$. Furthermore, \eqref{lie113} and
\eqref{rLiedef} show that $\nu: \mathfrak{h}_{r} \to
\overline{\mathfrak{h}}$ is also a 3-Lie algebra map which
finishes the proof.
\end{proof}

In order to provide the classification of all complements we
introduce the following:

\begin{definition}\label{equivLie}
Let $(\mathfrak{g}, \mathfrak{h}, \triangleright, \triangleleft,
\leftharpoonup, \rightharpoonup)$ be a matched pair of 3-Lie
algebras. Two deformation maps $r$, $R: \mathfrak{h} \to
\mathfrak{g}$ are called \emph{equivalent} and we denote this by
$r \sim R$ if there exists $\sigma: \mathfrak{h} \to \mathfrak{h}$
a $k$-linear automorphism of $\mathfrak{h}$ such that for any $x$,
$y\in \mathfrak{h}$:
\begin{eqnarray*}\label{equivLiemaps}
&&\sigma[u_1, u_2, u_3]-[\sigma(u_1), \sigma(u_2), \sigma(u_3)]\\
\notag&=&(\sigma(u_1), \sigma(u_2))\ppl R(\sigma(u_3))+c.p.+\sigma(u_1)\trl(R(\sigma(u_2)), R(\sigma(u_3)))+c.p.\\
&&-\sigma\Big((u_1, u_2)\ppl r(u_3)+c.p.+ u_1\trl(r(u_2), r(u_3))+c.p.\Big).
\end{eqnarray*}
\end{definition}

To conclude this section, the following result provides the answer
to the classifying complements problem  for 3-Lie algebras:

\begin{thm}\label{clasformelorLie}
Let $\mathfrak{g}$ be a 3-Lie subalgebra of ${E}$,
$\mathfrak{h}$ a $\mathfrak{g}$-complement of ${E}$ and
$(\mathfrak{g}, \mathfrak{h}, \triangleright, \triangleleft,
\leftharpoonup, \rightharpoonup)$ the associated canonical matched
pair. Then $\sim$ is an equivalence relation on the set $
{\mathcal D}{\mathcal M} \, ( \mathfrak{h}, \mathfrak{g} \, | \, (
\triangleright, \triangleleft, \leftharpoonup, \rightharpoonup )
)$. If we denote by ${\mathcal M}{\mathcal H}^{2} (\mathfrak{h}, \mathfrak{g}) \, := \, {\mathcal D}{\mathcal M} \,
(\mathfrak{h}, \mathfrak{g} \, | \, (\triangleright,
\triangleleft, \leftharpoonup, \rightharpoonup) )/\sim$, then we have
$$
{\mathcal M}{\mathcal H}^{2} (\mathfrak{h}, \mathfrak{g})
\longrightarrow {\mathcal F} (\mathfrak{g}, \, {E}),
$$
$$
\overline{r} \mapsto \mathfrak{h}_{r}
$$
is a bijection between ${\mathcal M}{\mathcal H}^{2}
(\mathfrak{h}, \mathfrak{g})$ and the
isomorphism classes of all $\mathfrak{g}$-complements of ${E}$.
In particular, the factorization index of $\mathfrak{g}$ in $E$ is
computed by the formula:
$$
[E: \mathfrak{g}] = | {\mathcal M}{\mathcal H}^{2}
(\mathfrak{h}, \mathfrak{g})|.
$$
\end{thm}

\begin{ex}
Let $E$ be the 6-dimensional 3-Lie algebra defined with respect to a basis
$\{x_1, x_2, x_3, u_1, u_2, u_3\}$ by the skew-symmetric bracket
\begin{eqnarray*}
&&[x_1,x_2,x_3] = x_1,\quad[u_1,u_2,x_2] = \alpha x_1+u_1,\\
&&[x_1,u_2,x_2] = u_1,\quad[u_1,x_2,x_3] = u_1,
\end{eqnarray*}
where $\alpha$ is an arbitrary parameter.

Let $\mathfrak{g}$ be the 3-Lie subalgebra of $E$ with basis $\{x_1, x_2, x_3\}$
and $\mathfrak{h}$ be the abelian 3-Lie algebra of dimension $3$ with basis $\{u_1, u_2, u_3\}$ . Then $\mathfrak{h}$ is a
$\mathfrak{g}$-complement of $E$ with the associated matched pair $(\mathfrak{g}, \, \mathfrak{h})$ given as follows:
$$
(u_1,u_2) \trr x_2= -\alpha x_1, \quad (u_1,u_2) \ppl x_2= u_1,\quad u_1\trl(x_2,x_3)= u_1,\quad u_2\trl(x_2,x_1)= u_1.
$$
It is easy to see that the map $r: \mathfrak{h} \to\mathfrak{g}$ given by
\begin{eqnarray*}
r(u_1) &=& 0, \quad r(u_2) =b_1x_1+ b_2x_2+b_3x_3,\\
r(u_3) &=& c_1x_1+ c_3x_3,
\end{eqnarray*}
is a deformation map  associated with the above matched pair of 3-Lie  algebras, where $b_i,c_i$ are arbitrary parameters.
Furthermore, the $r$-deformation of $\mathfrak{h}$ has the bracket given by
\begin{eqnarray*}
[u_1, u_2, u_3]_r=b_2c_3u_1.
\end{eqnarray*}
If $b_2c_3 \neq 0$, then  $\mathfrak{h}_{r}$ is not isomorphic to $\mathfrak{h}$ as $\mathfrak{h}_{r}$ is not abelian. Since we have only two types of 3-Lie algebras of dimension $3$, we obtain that $[E :
\mathfrak{g}] = 2$.
\end{ex}

\section{Conclusions and problems}
In this paper, the theory of extending structures and unified products for  3-Lie algebras are developed. We found that the  extending structures can be classified by using some non-abelian cohomology and deformation map.
There are many problems which are deserved to be considered in the future. Firstly, can all the results of this paper be generalized to the case of $n$-Lie algebras or  $n$-Leibniz algebras?
We believe that there is no essential difficulty beyond some complicated computations with respect to $n$-bracket.
Secondly,  if we consider  the special  matched pair $(\g,\g^*)$  where $\g^*$, the dual space of $\g$, is also a non-abelian 3-Lie algebra,  we will obtain the general theory of  "3-Lie bialgebras"
which is different from "3-Lie bialgebras" defined  in  \cite{Bai2,DBL}. The intrinsic relationship between the theory of unified products for  3-Lie algebras and 3-Lie bialgebras is deserved to be found.
Finally, an interesting problem is  to develop the theory of  flag extending structures for 3-Lie algebras as in \cite{AM1,AM3,AM4}.
The solutions of these problems are left to future investigations.

\subsection*{Acknowledgements}
We would like to thank the referee for careful reading and for valuable suggestions on this paper.
This work is supported in part by Natural Science Foundation of China (11501179, 11961049).
Part of this work was done while the author was visiting Courant Research Centre, Georg-August Universit\"{a}t G\"{o}ttingen from June to September, 2013. He is grateful to Professor Chenchang Zhu for invitation and hospitality.

\vskip7pt

\footnotesize{\noindent 
 College of Mathematics and Information Science, Henan Normal University, Xinxiang 453007, P. R. China;\\
 E-mail address:\texttt{{
 zhangtao@htu.edu.cn}}.

}

\end{document}